\documentclass[a4paper,leqno]{article}

\usepackage{RRA4}
\usepackage{a4wide}

\RRdate{Novembre 2010}

\RRauthor{
  Hend Ben Ameur%
  \thanks[1]{FSB and ENIT-LAMSIN, University of Tunis, BP 37, Le Belv\'ed\`ere,
    1002 Tunis, Tunisia.
    {\tt hbenameur@yahoo.ca}.}%
  \and Guy Chavent%
  \thanks[2]{Projet Estime.
    {\tt \{Guy.Chavent,Francois.Clement,Pierre.Weis\}@inria.fr}.}
  \and Fran\c{c}ois Cl\'ement%
  \thanksref{2}
  \and Pierre Weis%
  \thanksref{2}
}
\authorhead{H. Ben Ameur, G. Chavent, F. Cl\'ement, \& P. Weis}

\RRtitle{Segmentation d'image par indicateurs de raffinement multidimensionnels}
\RRetitle{Image Segmentation with Multidimensional Refinement Indicators}

\RRresume{

Nous transposons une technique de contr\^ole optimal \`a la segmentation
d'image.
L'id\'ee est de voir la segmentation d'image comme un probl\`eme d'estimation
de param\`etre o\`u le param\`etre \`a estimer est la couleur des pixels de
l'image.
Nous utilisons la technique de {\em param\'etrisation adaptative} qui construit
it\'erativement une repr\'esentation optimale du param\`etre en régions
uniformes formant une partition du domaine, et correspondant ainsi \`a une
segmentation de l'image.
Nous minimisons une fonction d'erreur au cours des it\'erations, et le
partitionnement de l'image en r\'egions est guid\'e de façon optimale par le
gradient de cette erreur.
L'algorithme de segmentation r\'esultant h\'erite des bonnes propri\'et\'es de
son origine contr\^ole optimal~: fondement, robustesse et flexibilit\'e.
}
\RRabstract{

We transpose an optimal control technique to the image segmentation problem.
The idea is to consider image segmentation as a parameter estimation problem.
The parameter to estimate is the color of the pixels of the image.
We use the {\em adaptive parameterization} technique which builds iteratively
an optimal representation of the parameter into uniform regions that form a
partition of the domain, hence corresponding to a segmentation of the image.
We minimize an error function during the iterations, and the partition of the
image into regions is optimally driven by the gradient of this error.
The resulting segmentation algorithm inherits desirable properties from its
optimal control origin: soundness, robustness, and flexibility.
}

\RRmotcle{segmentation d'image,
paramétrisation adaptative,
problème inverse,
contrôle optimal
}
\RRkeyword{image segmentation,
adaptive parameterization,
inverse problem,
optimal control
}

\RRprojet{Estime}

\RRdomaine{5}
\RRtheme{Observation et mod\'elisation pour les sciences de l'environnement}

\RCParis

\newcommand{\card}{{\rm card}}
\newcommand{\calF}{{\mathcal F}}
\newcommand{\calM}{{\mathcal M}}
\newcommand{\calP}{{\mathcal P}}

\newcommand{\bfR}{\mathbf{R}}
\newcommand{\egaldef}{\stackrel{\mathrm{def}}{=}}

\newcommand{\demi}{{\frac{1}{2}}}

\newcommand{\Cuttingstrategies}{Cutting strategies}
\newcommand{\cuttingstrategies}{cutting strategies}
\newcommand{\cuttingstrategy}{cutting strategy}
\newcommand{\Bestinafamily}{{\em Best in a family}}
\newcommand{\bestinafamily}{{\em best in a family}}
\newcommand{\Overallbest}{{\em Overall best}}
\newcommand{\overallbest}{{\em overall best}}

\newcommand{\multiscalarstrategies}{multiscalar strategies}
\newcommand{\multiscalarstrategy}{multiscalar strategy}
\newcommand{\Bestcomponentonly}{{\em Best component only}}
\newcommand{\bestcomponentonly}{{\em best component only}}
\newcommand{\Bestcomponentforeach}{{\em Best component for each}}
\newcommand{\bestcomponentforeach}{{\em best component for each}}
\newcommand{\Combinebestcomponents}{{\em Combine best components}}
\newcommand{\combinebestcomponents}{{\em combine best components}}

\newcommand{\quartps}{0.19}

\newtheorem{remark}{Remark}

\begin{document}

\RRNo{7446}
\makeRR


\section{Introduction}
\label{sec:intro}

The segmentation of an image is usually defined as the identification of
subsets of points with similar image properties, such as the brightness or
the color (local properties), or the texture (global property).
Image segmentation is usually subdivided into two main classes:
edge detection approaches that follow the variations of the image properties,
and region-based approaches that follow directly the image properties,
e.g. see~\cite{pav:spr:80,prw:sla:82,pp:ris:93,spi:sis:93,sk:cis:94,cd:ibf:95,fan:ais:01,bai:seg:03,mun:sis:03}.

Introducing an optimization process to guide a segmentation algorithm
has already been investigated.
In~\cite{klm:mai:94}, authors introduce an energy to minimize and develop a
``region merging'' algorithm for grey level and texture segmentation.
The energy method in image segmentation is introduced by Geman and
Geman~\cite{gg:ieee:84}.
Later a large variety of energy methods has been developed.
A formulation of the image segmentation as a piecewise modeling problem
combined with an energy method is given in~\cite{pfg:nce:06}.
Uchiyama and Arbib~\cite{ua:cis:94} incorporate optimization in a
merge-split, region-based image segmentation method; their algorithm
converges to a local optimum solution of clustering color space based on the
least sum of squares criterion.

The connection between inverse problems and image segmentation has already
been noticed, e.g. see~\cite{bon:ipi:08}.
Furthermore, an approach close to the one we are developing can be found
in~\cite{wu:asm:93}, where the segmentation is built through a piecewise
least-squares approximation of a function related to image intensity.
But unlike ours, the refinement of the segmentation is uniform and the
refinements are selected from a fixed predefined set.

The {\em optimal adaptive segmentation algorithm} presented in this paper is
also based on an optimization process.
At each step, the refinement is optimally obtained from the gradient of the
objective function to be minimized.
The resulting regions are not geometrically constrained, as they are not
required to be connected.
Noticeably, all the steps of the algorithm are meaningful, since the
algorithm progressively goes from the coarsest segmentation to the finest
one, i.e. from a single monochrome region up to the complete reconstruction
of the original image with one region per distinct color.

This leads to some kind of {\em segmentation on demand}: since the $n$-th
iteration provides the segmentation into $n$ regions, the algorithm can
segment the image into any given number of zones, with each zone having an
optimal color with respect to the overall error from the original image.
Moreover, since each step of the algorithm selects a single region
and splits it into two parts, the adaptive segmentation algorithm can be used
to {\em interactively} refine the segmentation of an image.

This adaptive approach was initially developed for the inverse problem of
parameter estimation: given a forward modeling map, the
{\em multidimensional adaptive parameterization} algorithm~\cite{ben:mri:08}
solves the inverse problem of estimating a multidimensional distributed
parameter from measured data.
For instance, colors are defined by their three Red/Green/Blue (RGB)
components, thus color images can be seen as three-dimensional distributed
parameters.
Hence, when the parameter to be identified is a color image and the forward
map to be inverted is the identity, this optimal control technique provides
naturally an image segmentation algorithm that refines regions guided by the
misfit between the original and segmented image.

As illustrated in~\cite{ben:mri:08}, brute force application of this
algorithm to image segmentation gives already interesting results.
But the original handling of the multidimensional aspect of the parameter
proposed there was rather naive.
So this aspect is revisited, developed and enriched in the present paper
specifically for the case of image segmentation:
beside {\em vector segmentation} algorithms, where each region is associated
with one color, we also introduce {\em multiscalar segmentation} algorithms,
where each color channel has its own segmentation.
When transposed back to the general case of nonlinear inverse problem, the
new algorithms are also expected to provide more flexibility.

Notice that although adaptive parameterization uses so-called {\em cuttings}
that can be considered as defining edges, this ``edge detection'' aspect is
just an auxiliary step to define the partitioning of the set of pixels into
regions.
Thus, we still classify our image segmentation algorithm in the region-based
techniques.

\bigskip

The paper is organized as follows.
In Section~\ref{sec:segmentation}, we explain the similarity between
parameter estimation and image segmentation and we illustrate the usage of
adaptive parameterization to segment a simple image.
In Section~\ref{optimal_adaptive_segmentation}, we give some mathematical
foundations and the description of the optimal adaptive {\em vector} and
{\em multiscalar} segmentation algorithms.
In Section~\ref{s:numerical_results}, we use the previous algorithms to
segment a color image, and compare their relative performances.


\section{From parameter estimation to image segmentation}
\label{sec:segmentation}

\subsection{Parameter estimation}
\label{parameter_estimation}

Parameter estimation is a class of inverse problems, where a simulation model
has been constructed but some of its parameters are not precisely known.
It consists in trying to recover these parameters from a collection of
input-output experimental data.

A classical approach to this inverse problem is to formulate it as the
minimization of an objective function defined as a least-squares misfit between
experimental output data and the corresponding model-generated quantities:
\begin{equation}
  \label{eqn:lsp}
  \mbox{minimize } J (c) = \demi \| d - \calF (c) \|^2,
\end{equation}
where~$c$ denotes the unknown parameter, $d$~the data and~$\calF$ is an
operator corresponding to the model describing the cause-to-effect link
from~$c$ to~$d$.

For a large class of inverse problems, the unknown parameters are space
dependent coefficients in a system of partial differential equations modeling a
physical problem.
Due to the difficulty and high cost of experimental measurements, the available
data are often too scarce to allow for the determination of the parameter value
in each cell of the computational grid.

Therefore, one has to reduce the number of unknown parameters:
a convenient way is to search for unknown distributed parameters as piecewise
constant functions, as an approximation to the necessarily more complex but
unattainable reality.

Then parameter estimation amounts to identify both the spatial distribution of
{\em regions} of constant value and the associated {\em parameter values}.
The result of the parameter estimation is then a {\em partition} of the spatial
domain into regions together with
{\em one (possibly multi-dimensional) parameter value per region}.

\subsection{Image segmentation}
\label{image_segmentation}

Let us now consider the case where the data~$d$ is an image, made of a set of
pixels (picture elements).
Each pixel is the association of a cell~$x_i$ with its corresponding color
value~$d_i$:
\begin{equation}
  \label{eqn:image}
  d = \{ (x_i, d_i), i \in I \},
  \quad \mbox{where~$I$ is the set indexing the image pixels}.
\end{equation}
We shall denote by~$D$ the domain of the image, i.e. the set of all its
cells~$x_i$:
\begin{equation}
  \label{eqn:domain}
  D = \{ x_i, i \in I \}.
\end{equation}

\paragraph*{Vector segmentation.}

A {\em vector segmentation} of~$d$ consists in replacing it by an image
$c=\{(x_i,c_i),i{\in}I\}$ (where~$c$ stands for ``computed'') with the
same domain such that:
\begin{itemize}
\item $c$~has a uniform color~$s_\bfR$ (where~$s$ stands for ``segmented'')
  on each region~$\bfR$ of a partition~$\calP$ of its domain~$D$:
  \begin{equation}
    \label{eqn:segmented_c}
    c_i = s_\bfR \mbox{ as soon as } x_i \in \bfR,
  \end{equation}
\item the image~$c$ bears some resemblance to the original image~$d$.
\end{itemize}
Image {\em (vector) segmentation} methods are categorized into three main
classes: pixel based, edge based and region based
approaches~\cite{che:cis:01,spi:sis:93}.
Edge detection is mainly used for gray level image segmentation, it consists in
locating points with abrupt changes in gray level, in~\cite{nader:nsh:08}
authors give a survey on techniques for edge detection.
Region based approaches include region growing, region splitting, region
merging and their combination.
They are also divided into two classes: top-down and
bottom-up~\cite{bor:bu:08}.
A combination of region growing and region merging techniques is presented
in~\cite{tre:tb:97}.
In~\cite{ohta:oks:80,tom:86}, authors use a homogeneous criterion based on a
one-dimensional histogram to develop a region splitting technique.

Here we have chosen to ensure the resemblance between the segmented and
original images~$c$ and~$d$ by requiring that their least-squares misfit is
made ``small'', or, in more mathematical terms:
\begin{equation}
  \label{eqn:J_c}
  \mbox{the segmented image~$c$ minimizes }
  J (c) = \demi \sum_{i \in I} \| d_i - c_i\|^2
\end{equation}
under the constraint that the number of regions of~$\calP$ is
``small compared to the number of pixels of~$d$''.

Without the constraint on the number of regions, a trivial exact solution to
problem~(\ref{eqn:J_c}) is the segmentation made of one region for each pixel,
with the region color being the color of its unique pixel.
However, this trivial exact solution is not satisfactory, since it is not
adequate to the segmentation idea.
It corresponds to the ultimate level of {\em over-segmentation}.

This approach can be seen as a top-down region based one since we keep refining
the regions: at each iteration, the algorithm detects a discontinuity in one
region color $s_\bfR$ which give a smaller value to $J(c)$, so the algorithm
can also be seen as defining an implicit edge detection.

The result of a {\em vector segmentation} of an image~$d$ is hence a
{\em partition}~$\calP$ of its domain~$D$ into regions, together with the data
of {\em one color value~$s_j$ per region}.
In such a segmentation, all color components use the same partition.

\paragraph*{Multiscalar segmentation.}

Another way to obtain an image with piecewise constant colors is to find a
different partition for each color component.
We have chosen to use the RGB basis, but any other representation could be used
as well.
So we define a {\em multiscalar segmentation} of~$d$ as the operation which
replaces~$d$ by an image $c=\{(x_i,c_i),i{\in}I\}$ such that:
\begin{itemize}
\item for each $k=R,G,B$, the component~$c^k$ of~$c$ has a uniform
  intensity~$s_\bfR^k$ on each region~$\bfR$ of a partition~$\calP^k$ of the
  image domain~$D$:
  \begin{equation}
    \label{eqn:segmented_ck}
    \mbox{for } k = R, G, B
    \mbox{ one has: } c_i^k = s_\bfR^k
    \mbox{ as soon as } x_i \in \bfR,
  \end{equation}
\item the image~$c$ bears some resemblance to the original image~$d$,
\end{itemize}
in the sense that~$c$ solves the same optimization problem:
\begin{equation}
  \label{eqn:Jk_c}
  \mbox{the segmented image~$c$ minimizes }
  J (c) = \sum_{k = R, G, B} J^k (c)
\end{equation}
where $J^k(c)=\demi\sum_{i{\in}I}(d_i^k-c_i^k)^2$, under the same constraint
that the number of regions of the partitions~$\calP^R$, $\calP^G$
and~$\calP^B$ is ``small compared to the number of pixels of~$d$''.

\bigskip

The approach proposed here can produce both {\em vector} and {\em multiscalar}
segmentations, and controls the segmentation level.

\subsection{Image segmentation considered as a particular
  parameter estimation problem}
\label{segmentation_as_estimation}

Similarities between the two previous sections show that vector segmentation of
images can be seen as a parameter estimation problem, where the forward
model~$\calF$ reduces to the identity map (compare~(\ref{eqn:lsp})
and~(\ref{eqn:J_c})) and complete data are available (since the exact color is
known for each cell).

The parameter estimation problem, because of the nonlinearity of the governing
model~$\calF$ and the scarcity of data, is often ill-posed and difficult to
solve.
In contrast, one can expect that the parameter estimation algorithms will work
efficiently for vector segmentation, because of the good properties of the
associated forward map~$\calF$ and data.
So we adapt in this paper the multidimensional refinement indicators
algorithm~\cite{ben:mri:08} to the specificities of vector segmentation, and
extend it to the determination of multiscalar segmentations.

\begin{remark}
  \label{rem:texture}
  The segmentation of an image can also be attempted according to other
  properties than color, for example texture.
  The approach developed in this paper can be adapted to these situations by a
  proper choice of the objective function~(\ref{eqn:J_c}) or~(\ref{eqn:Jk_c}).
\end{remark}

\subsection{Illustration on an example}
\label{pedagogical_example}

In the following, we apply the adaptive parameterization algorithm
of~\cite{ben:mri:08} to perform vector segmentation of simple images.

\begin{figure}[ht]
  \begin{center}
    (a)
    \includegraphics[width=\quartps\textwidth]{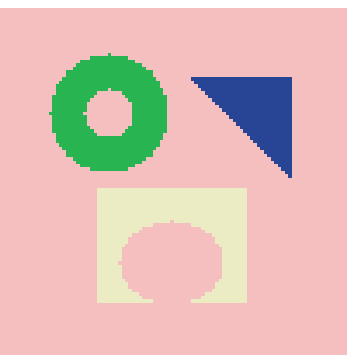} \\[0.2cm]
    (b)
    \includegraphics[width=\quartps\textwidth]{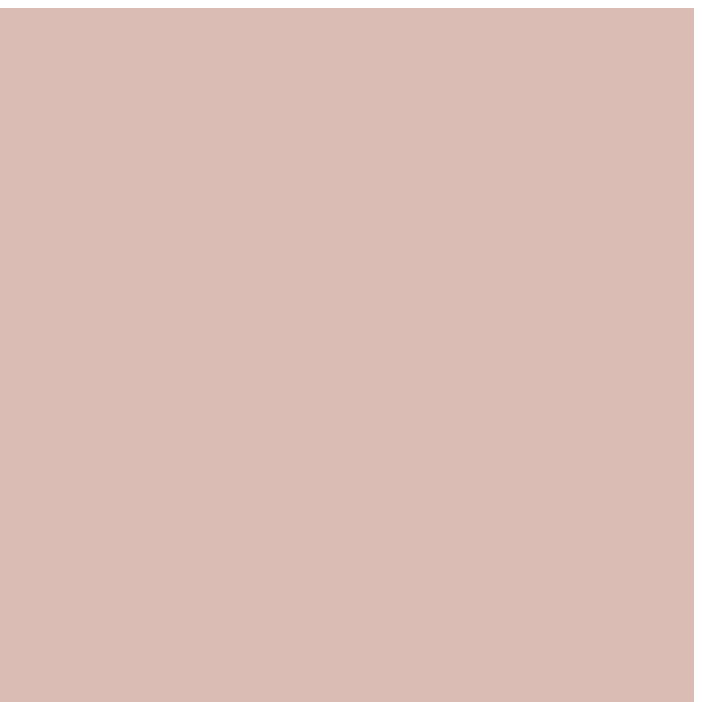}
    (c)
    \includegraphics[width=\quartps\textwidth]{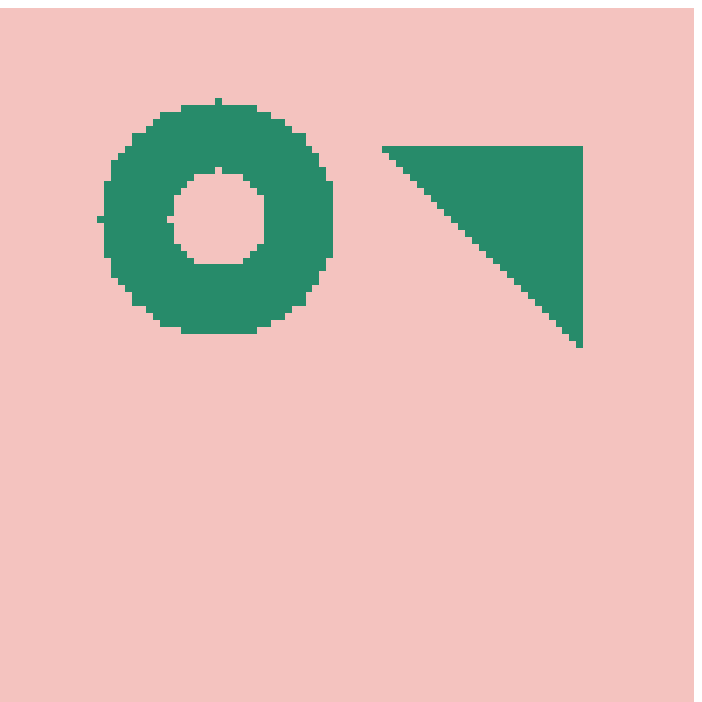}
    (d)
    \includegraphics[width=\quartps\textwidth]{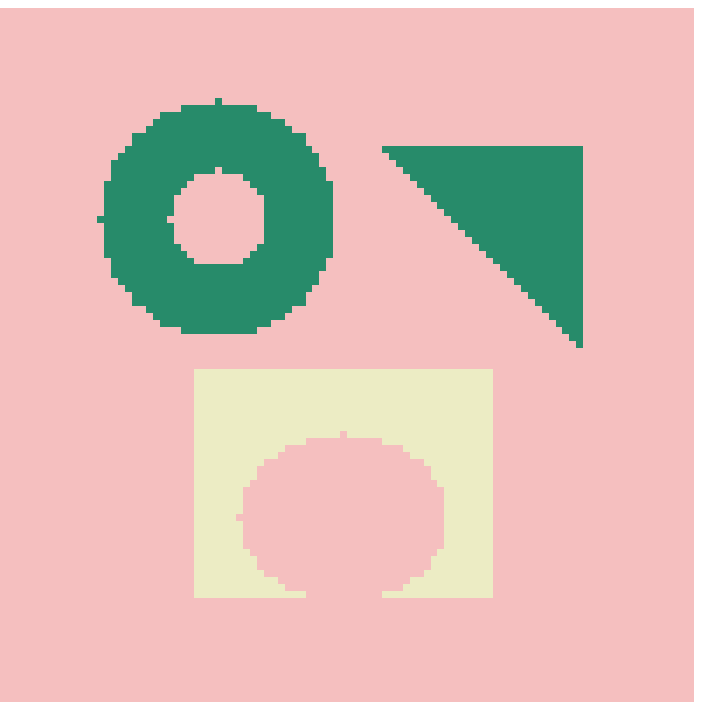}
    (e)
    \includegraphics[width=\quartps\textwidth]{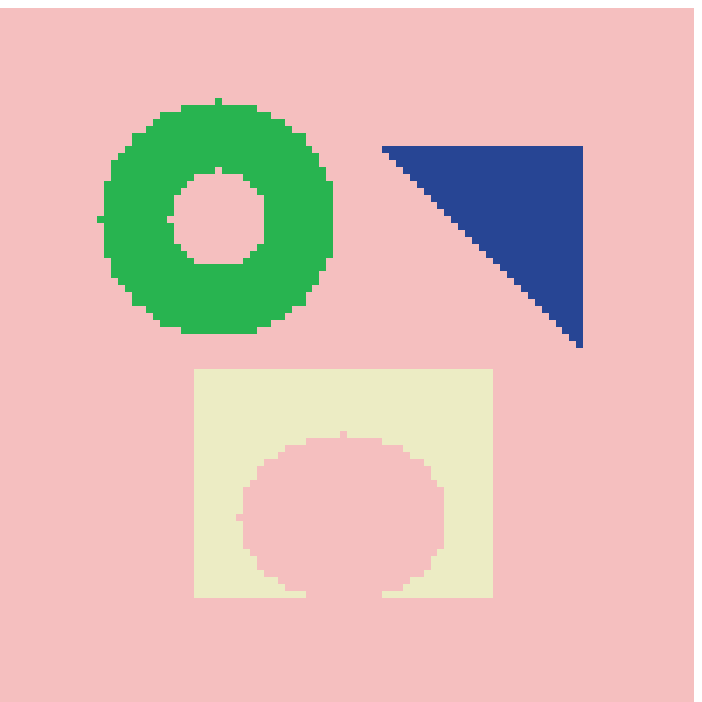}
    \caption{A simple example:
      (a)~is the image to segment, and~(b), (c), (d), and~(e) are the four
      successive iterations of the adaptive parameterization algorithm.}
    \label{fig:simple}
  \end{center}
\end{figure}

In the first example, Figure~\ref{fig:simple}(a) is the original image~$d$,
i.e. the {\em data} of the parameter identification problem.
It is made of three geometrical patterns with color green, yellow and blue on
top of a pink background.
During the iterations illustrated from Figure~\ref{fig:simple}(b)
to Figure~\ref{fig:simple}(e), we build the segmentation in an iterative way.
We start with one (homogeneous) region in Figure~\ref{fig:simple}(b),
it is split into two regions in Figure~\ref{fig:simple}(c),
a new region is detected in Figure~\ref{fig:simple}(d),
and in Figure~\ref{fig:simple}(e) the algorithm has converged to a segmentation
made of four regions, in each region the vector color is constant.

This test case exemplifies the relevance of all the iterations.
All iterations are meaningful:
the first step gives a monochrome image with the optimal color (the mean color
of the picture elements);
the second step gives the best splitting of the image into two zones with the
best two colors;
and so on.

\begin{figure}[ht]
  \begin{center}
    (a)
    \includegraphics[width=\quartps\textwidth]{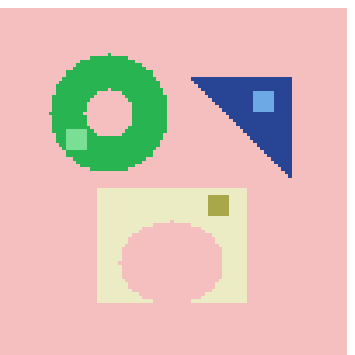} \\[0.2cm]
    (b)
    \includegraphics[width=\quartps\textwidth]{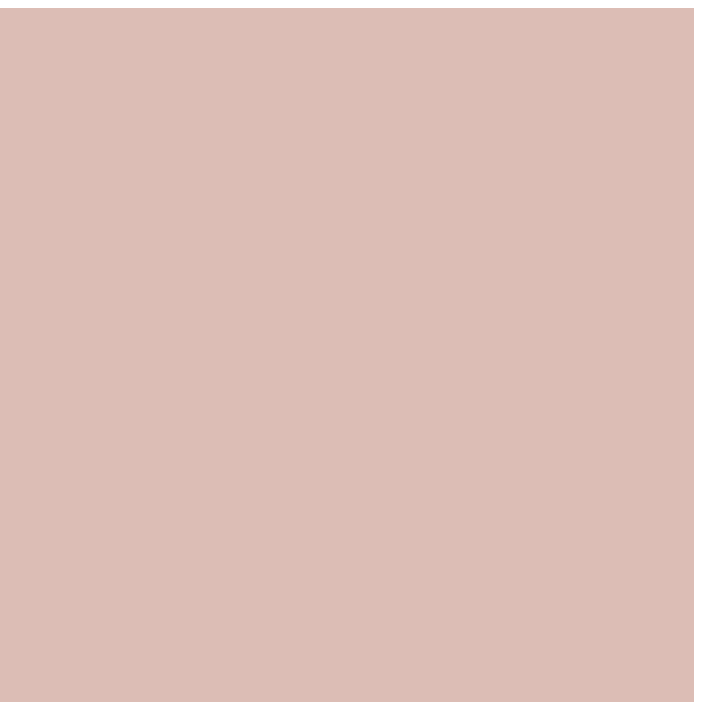}
    (c)
    \includegraphics[width=\quartps\textwidth]{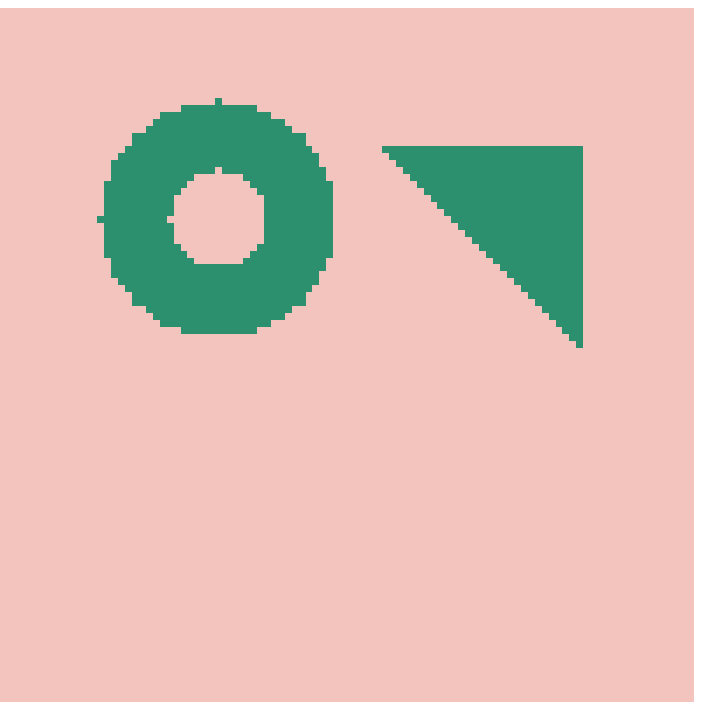}
    (d)
    \includegraphics[width=\quartps\textwidth]{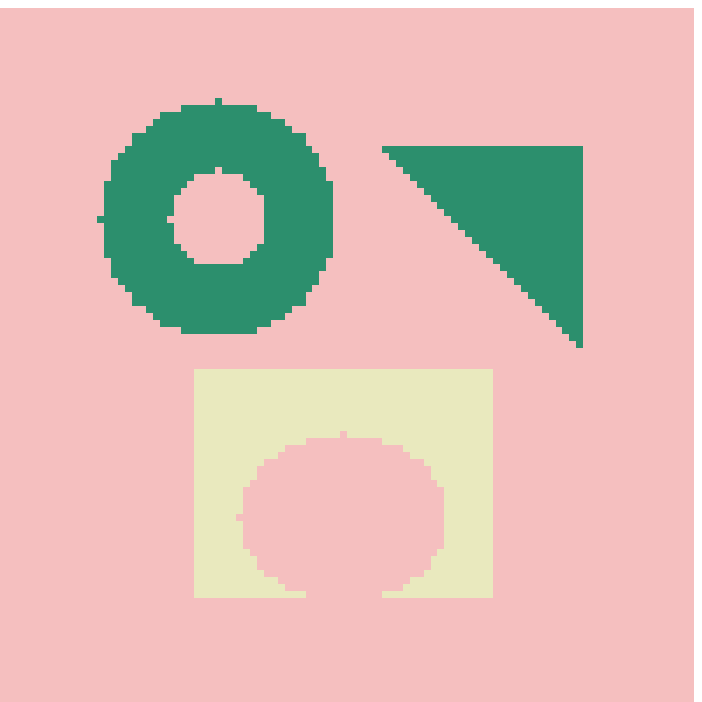}
    (e)
    \includegraphics[width=\quartps\textwidth]{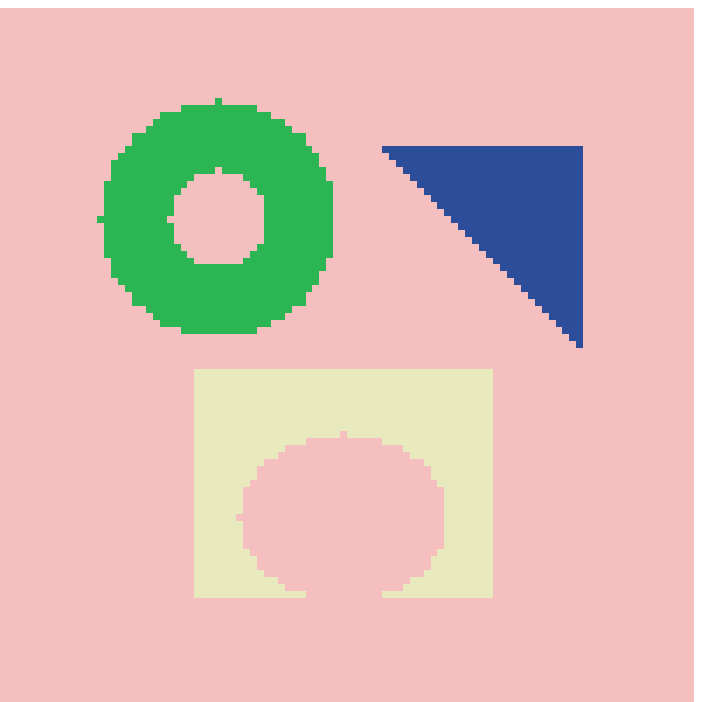} \\[0.2cm]
    (f)
    \includegraphics[width=\quartps\textwidth]{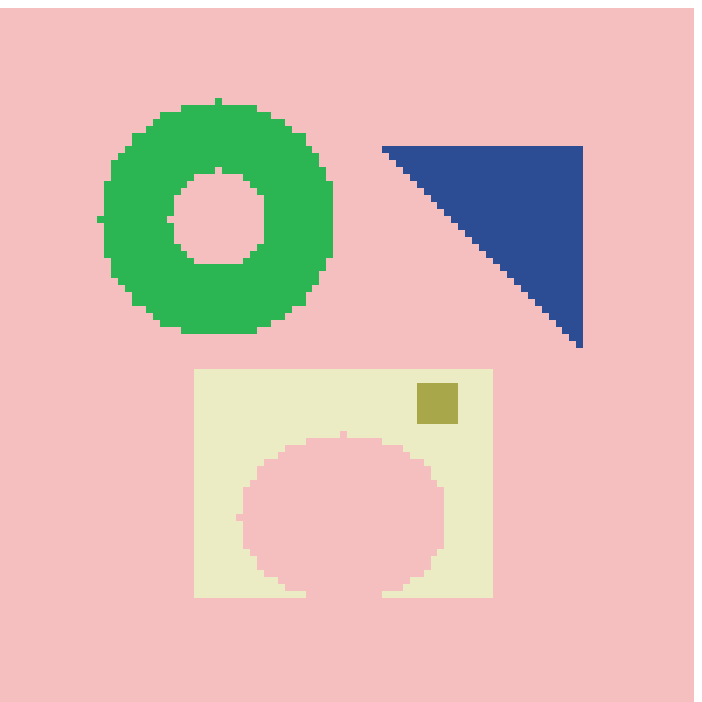}
    (g)
    \includegraphics[width=\quartps\textwidth]{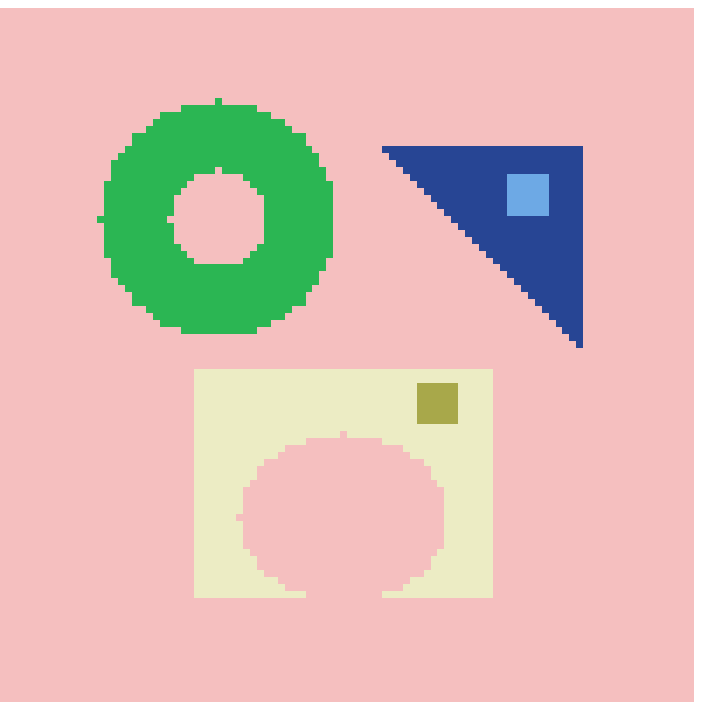}
    (h)
    \includegraphics[width=\quartps\textwidth]{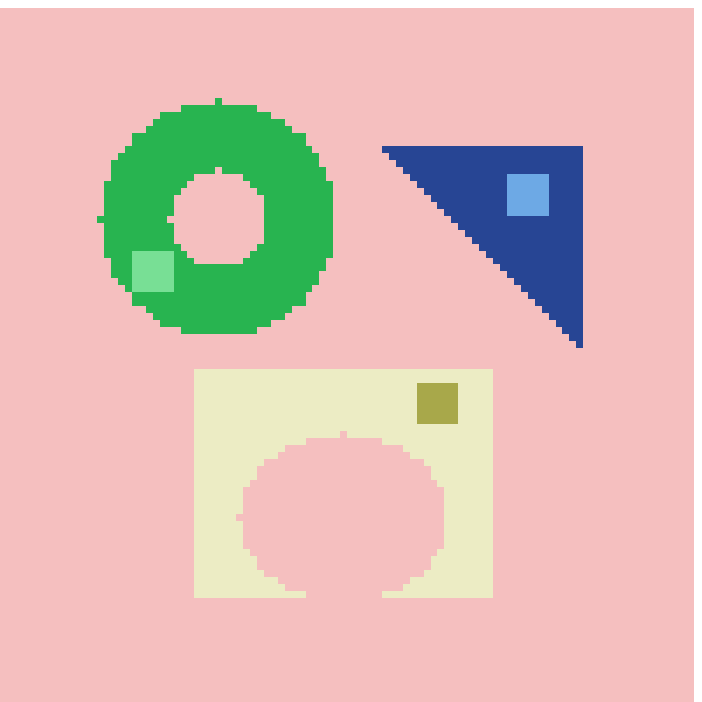}
    \caption{A perturbation of the simple example:
      (a)~is the image to segment, and~(b), (c), (d), (e), (f), (g), and~(h)
      are the seven successive iterations of the adaptive parameterization
      algorithm.}
    \label{fig:perturb}
  \end{center}
\end{figure}

In the second example, the image to be segmented Figure~\ref{fig:perturb}(a)
is a perturbation of the previous image Figure~\ref{fig:simple}(a) with small
square inclusions in the yellow, blue, and green patterns.
Iteration results are illustrated from Figure~\ref{fig:perturb}(b) to
Figure~\ref{fig:perturb}(h).
The four first iterations do not detect the perturbations, and the fourth
Figure~\ref{fig:perturb}(e) is exactly the final segmentation of the non
perturbed image (to be correct, only the partition associated with it is the
same).
The last three iterations successively detect the three perturbations.

This second example shows that the larger the contrast is, the faster it is
detected, which means that a stopping criterion could be the number of regions
supposed to be the more relevant.
Also notice that the data image is retrieved in only seven iterations, which is
exactly the number of distinct colors in the perturbed data image.


\section{Optimal adaptive segmentation}
\label{optimal_adaptive_segmentation}

{\em Adaptive segmentation} algorithms have been proposed by many authors.
In~\cite{ bbb:ras:00} level sets are used to design a robust evolution model
based on adaptation of parameters with respect to data.
In~\cite{rck:ass:00} the adaptivity is related to the fact that the
segmentation is adapted to different regions:
regions are segmented using an appropriate method according to their types.
The adaptivity is achieved through the variation of the size of the
neighborhoods used for the estimation of the intensity functions in
\cite{pac:uas:01}.
Finally, adaptive segmentation in~\cite{che:ais:02} is related to the use of an
adaptive clustering algorithm~\cite{pap:aca:92} to obtain the specially
adaptive dominant colors.

The starting point for the method proposed in this paper is the adaptive
multidimensional parameter estimation algorithm proposed in~\cite{ben:mri:08}
for the general nonlinear inverse problem of
section~\ref{parameter_estimation}.
It proceeds iteratively, using multidimensional refinement indicators to split
at each step one of the regions where the multidimensional parameter is
constant.
This approach allows to control the issue of {\em over-parameterization} by
stopping the refinement process when the objective function is of the magnitude
of the uncertainty or noise level.

As indicated in section~\ref{segmentation_as_estimation}, we shall adapt this
algorithm to the specificities of image segmentation.
The result appears as an iterative region-based {\em vector segmentation}
method, initialized by a one-region segmented image of uniform color.
At each iteration, it adapts the segmentation by splitting one region.
The region to be split and the location of the splitting, taken from a set of
user-defined admissible cuttings, are chosen according to values of {\em exact}
or {\em first order} {\em refinement indicators} as defined
in~\cite{cha:nllsip10,ben:mri:08}.
These indicators are computed from the objective function~(\ref{eqn:J_c})
or~(\ref{eqn:Jk_c}) or its gradients with respect to color components.
Once the regions are given, the {\em color} in each region is computed by a
straightforward minimization process.
With this choice, the new segmented image is likely to produce a significant
decrease of the misfit function.
The algorithm is stopped when the desired number of regions is attained,
i.e. when the segmentation is fine enough to retain the information of
interest, and coarse enough to avoid over-segmentation.

Note that regions need not to be necessarily connected in our approach.
This feature may be considered a major drawback for some image segmentation
applications that requires connected regions only.
Extending our algorithm to deal with this aspect is an interesting research
problem that has still to be investigated.

We describe in the next sections the specification of the adaptive
parameterization technique to the {\em vector segmentation}
problems~(\ref{eqn:J_c}), and its extension to the
{\em multiscalar segmentation} problem~(\ref{eqn:Jk_c}).
We recall that the original image is denoted by $d=\{(x_i,d_i),i{\in}I\}$
(data), the segmented image by $c=\{(x_i,c_i),i{\in}I\}$ (computed).

\subsection{Case of vector segmentation}

Let~$\calP$ be the current partition.
The associated segmented image~$c$, solution of the least-squares
problem~(\ref{eqn:J_c}) for the partition~$\calP$, is given by
$c=\calM_\calP\tilde{\calM}_\calP\,d$, where
$\calM_\calP:\{s_\bfR,\bfR\in\calP\}{\mapsto}c$ defined
by~(\ref{eqn:segmented_c}) denotes the mapping that paints the pixels of each
region~$\bfR$ in~$\calP$ with a given color~$s_\bfR$, and
$\tilde{\calM}_\calP=(\calM_\calP^T\calM_\calP)^{-1}\calM_\calP^T$ is its
least-squares pseudoinverse, which computes the mean value in each region.
Hence once the partition~$\calP$ has been chosen, the best segmented image~$c$
is simply obtained by replacing on each region~$\bfR$ the image color by its
mean value.

So we concentrate now on the determination of a partition~$\calP$ with a
controlled number of regions and such that~$J(c)$ is as small as possible.

\subsubsection{Refinement indicators: how good is a cutting?}
\label{sss:ri_vector}

Let us first denote by~$C$ the cutting which splits a region~$\bfR$ of~$\calP$
into two subregions~$\bfR_+$ and~$\bfR_-$.
In order to evaluate the effect of this additional degree of freedom on the
decrease of the objective function~$J(c)$ in the vector segmentation
problem~(\ref{eqn:J_c}), one can use refinement
indicators~\cite[pp.~110--111]{cha:nllsip10}:
\begin{description}
\item[Exact indicators.]
  They consist in
  \begin{equation}
    \label{eqn:exact_indicators_norm}
    \Delta J = J (c) - J (c_C),
  \end{equation}
  where~$c_C$ is the solution of problem~(\ref{eqn:J_c}) for the
  partition~$\calP_C$ obtained by splitting region~$\bfR$ of~$\calP$
  into~$\bfR_+$ and~$\bfR_-$.
  These indicators are called ``nonlinear'' in~\cite{cha:nllsip10} because in
  the general case where~$\calF$ is nonlinear, the evaluation of~$c_C$
  requires the resolution of the full nonlinear least-squares
  problem~(\ref{eqn:lsp}) with the tentative refined partition~$\calP_C$.
  This prohibits their use for the test of a large number of cuttings.
  
  On the contrary, for our problem where~$\calF$ reduces to the identity map,
  $J$~is additive with respect to regions of~$\calP$, so that~$c_C$
  {\em coincides with~$c$ outside of}~$\bfR$.
  Its value on~$\bfR$ is given by the simple formulas:
  \begin{equation}
    \label{eqn:cCi}
    \left( c_C \right)_i =
    \left\{
      \begin{array}{lll}
        c_+ & \mbox{for } x_i \in \bfR_+ &
        \mbox{where~$c_+$ = mean value of } \{ d_j, x_j \in \bfR_+ \},\\
        c_- & \mbox{for } x_i \in \bfR_- &
        \mbox{where~$c_-$ = mean value of } \{ d_j, x_j \in \bfR_- \}.
      \end{array}
    \right.
  \end{equation}
  Hence~$\Delta J$ can be computed at a very low cost.
  This makes it possible to use exact indicators to rank a quite large number
  of tentative cuttings.
  
\item[First order indicators,] as defined in~\cite{ben:mri:08} for
  multidimensional parameters, are given by:
  \begin{equation}
    \label{eqn:first_order_indicators_norm}
    \| \lambda \|_q, \quad 1 \leq q \leq \infty, \quad
    \lambda = (\lambda^R, \lambda^G, \lambda^B),
  \end{equation}
  where the~$q$ is user defined and
  \begin{equation}
    \label{eqn:first_order_indicators_RGB}
    \lambda^k =
    \sum_{x_i \in \bfR_+} \frac{\partial J^k}{\partial c^k_i} -
    \sum_{x_i \in \bfR_-} \frac{\partial J^k}{\partial c^k_i}
    \quad \mbox{for } k = R, G, B.
  \end{equation}
  Because of the simple form of~$J^k$, the following relation holds for
  each color component (see~\cite{ben:mri:08}, or~\cite[p.~126]{cha:nllsip10}
  for the proof):
  \begin{equation}
    \label{eqn:identite_remarquable}
    \Delta J^k =
    \frac{(\lambda ^k)^2}{8} \, \frac{p}{p_+^k p_-^k}
    \quad \mbox{for } k = R, G, B,
  \end{equation}
  where~$p$ is the number of points~$x_i$ in~$\bfR$ and~$p_+^k$ and~$p_-^k$
  are the numbers of points~$x_i$ in the tentative subregions~$\bfR_+$
  and~$\bfR_-$.
  
  \begin{figure}[htb]
    \begin{center}
      \includegraphics[width=0.8\textwidth]{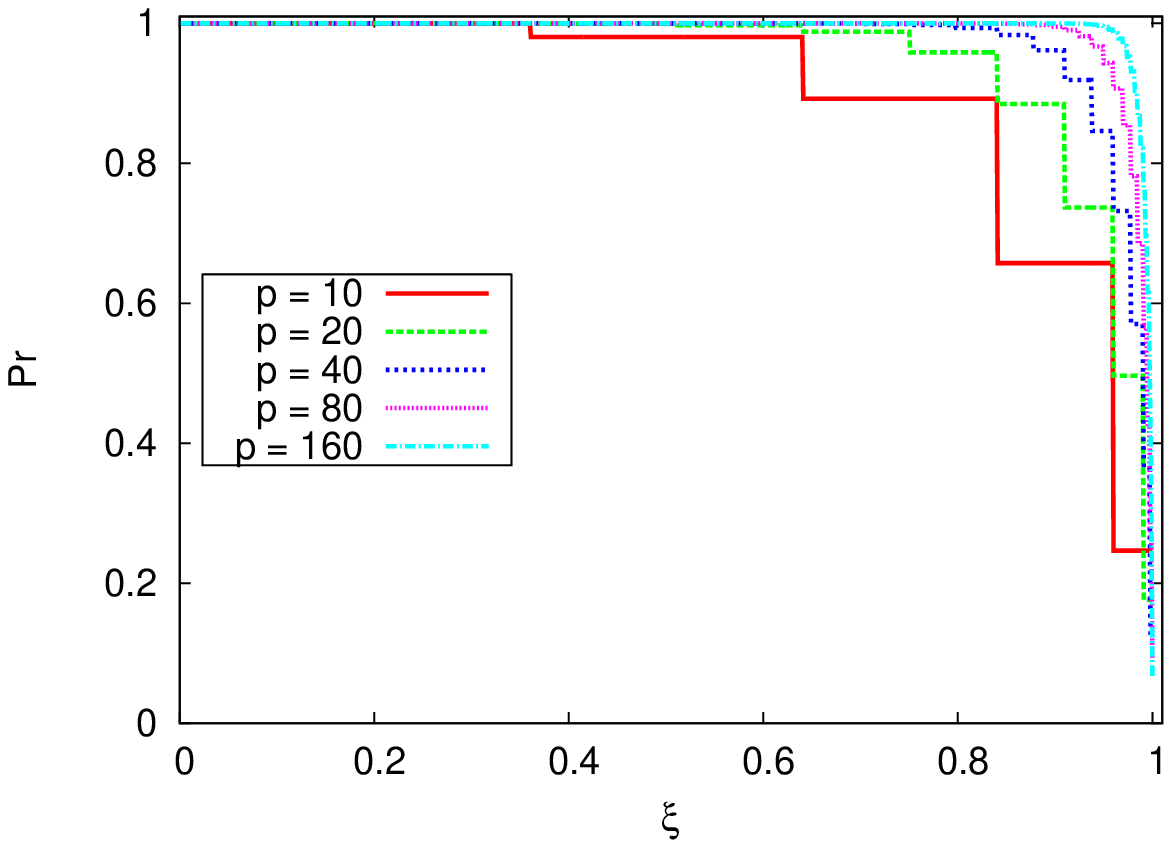}
      \caption{Comparison of exact and first order indicators for various~$p$.
        Each curve corresponds to the $\xi{\mapsto}Pr$ function defined
        in~(\ref{eqn:tau_to_pr}).}
      \label{fig:pr}
    \end{center}
  \end{figure}
  
  In order to illustrate the quality of this first order indicator, we have
  plotted in Figure~\ref{fig:pr} the curves%
  \footnote{%
    Since in~(\ref{eqn:identite_remarquable}), we have $p=p_+^k+p_-^k$,
    we also have
    $\frac{p_+^kp_-^k}{p^2}=
    \frac{p_+^k}{p}\left(1-\frac{p_+^k}{p}\right)\leq\frac{1}{4}$,
    and thus ${\Delta J}^k\geq\frac{(\lambda ^k)^2}{2p}$.
  }:
  \begin{equation}
    \label{eqn:tau_to_pr}
    \xi \in [0, 1] \mapsto Pr \left\{
      \xi \leq \frac{(\lambda ^k)^2}{2p} \Big / \Delta J^k \leq 1
    \right\}
  \end{equation}
  for regions~$\bfR$ with different numbers~$p$ of pixel, under the hypothesis
  that all partitions of~$\bfR$ are equiprobable.
  As one can see, this probability is practically one for~$\xi=0.9$ as soon
  as~$\bfR$ has more than 80~pixels.
  This shows that the first order indicators~$|\lambda^k|$ are well suited to
  rank cuttings inside a given region~$\bfR$.
  It also shows that they are not suited to rank cuttings belonging to
  different regions, unless they are replaced by $|\lambda^k|/p^\demi$,
  where~$p$ is the number of pixel of the region subject to the cutting.
\end{description}

The computational cost of first order indicator is here similar to that of
exact indicators, so one can wonder why to use an approximation when an exact
calculation is available at the same price.
The answer will appear in the next section.

\subsubsection{{\Cuttingstrategies}}
\label{sss:cutting_strategies}

Now that we know how to rank cuttings, the next step is to specify the
{\em {\cuttingstrategy}}, i.e. the family~$T$ of cuttings~$C$ to be tested for
the splitting of each region.

For nonlinear problems, adaptive parameterization has been shown to be
regularizing when the family~$T$ is adapted to the problem~\cite{benkal}.
But if enriching the cutting family always provides a better local decrease of
the cost function, it can also decrease the regularizing effect by leading the
algorithm into a local minimum.

However, for our segmentation problem, the forward map~$\calF$ is the identity,
so there is no danger of local minima, and we are free to choose the family~$T$
at our convenience, even if it is extremely large.
So we shall consider two {\em {\cuttingstrategies}}:
\begin{description}
\item[{\Bestinafamily}.]
  The set~$T$ of cuttings used for the tentative splitting of a region is a
  {\em user-defined family} of cuttings.
  It may incorporate any a priori information on the parameter when such
  information is available.
  
  We have used in our numerical results the set~$T$ made of all {\em vertical}
  and {\em horizontal} cuttings, which split any region~$\bfR$ in two
  subregions separated by one vertical or horizontal boundary.
  This family of cuttings has proven useful for the nonlinear parameter
  estimation problem, where one seeks to describe the parameter by regions with
  simple forms (see for example~\cite{ben:rci:02}, where~$T$ contained also
  diagonal and checkerboard cuttings).
  But we shall see in the numerical results section that it seems much less
  adapted to image segmentation.
  
  In any case, the number of cuttings to be tested with this strategy is
  bounded from above by the sum of the number of horizontal and vertical pixels
  in the image to be segmented.
  Hence the best tentative cutting~$C_\bfR^\star$ of region~$\bfR$ can be
  determined by an exhaustive computation of all
  {\em exact indicators}~${\Delta}J$ associated with all cuttings~$C$ of~$\bfR$
  in the chosen family~$T$, as they are not more expensive as the first order
  ones:
  \begin{equation}
    \label{eqn:best_cutting}
    \mbox{$C_\bfR^\star$ maximizes~${\Delta}J$ over all cuttings~$C{\in}T$
      of~$\bfR$}.
  \end{equation}
  
\item[{\Overallbest}.]
  Here~$T$ is made of {\em all possible 2-partitions} of the region~$\bfR$.
  So the number of cuttings to be tested for the determination
  of~$C_\bfR^\star$ is~$2^{\card(\bfR)}$, which becomes extremely large as soon
  as~$\bfR$ contains enough points.
  Hence one cannot afford to perform an exhaustive computation of all
  corresponding indicators, let them be exact or first order.
  So this {\cuttingstrategy} can be used only in cases where an analytical
  determination of the cut~$C_\bfR^\star$ which maximizes~${\Delta}J$
  or~$\|\lambda\|_q$ over~$T$ is available.
  To our best knowledge, such analytical solution is available neither for the
  exact indicator~${\Delta}J$, nor for the refinement
  indicators~$\|\lambda\|_q$ for any~$q$, except for~$q=\infty$.
  This is where the first order indicators play a crucial role in image
  segmentation.
  
  We describe now the analytical determination of $C_\bfR^\star$
  when~$q=\infty$.
  As it was noticed in~\cite{ben:mri:08},
  formula~(\ref{eqn:first_order_indicators_RGB}) shows that, for each color
  component~$k$, the cutting~$C_\bfR^{\star,k}$ which
  maximizes~$|\lambda_\bfR^k|$ is obtained by splitting~$\bfR$ according to the
  sign of $\partial J^k/\partial c_i^k$:
  \begin{equation}
    \label{eqn:first_order_cutting_opt_RGB}
    C_{\bfR_+}^{\star, k} = \left\{
      x_i \in \bfR \Big / \frac{\partial J^k}{\partial c_i^k} \geq 0
    \right\}, \quad
    C_{\bfR_-}^{\star, k} = \left\{
      x_i \in \bfR \Big / \frac{\partial J^k}{\partial c_i^k} < 0
    \right\}.
  \end{equation}
  The associated first order indicators of largest modulus for the splitting
  of~$\bfR$ are then, for each component:
  \begin{equation}
    \label{eqn:first_order_indicators_opt_RGB}
    \lambda_\bfR^{\star, k} = \left| \lambda_\bfR^{\star, k} \right| =
    \sum_{x_i \in \bfR} \left| \frac{\partial J^k}{\partial c_i^k} \right|
    \quad \mbox{for } k = R, G, B.
  \end{equation}
  So if we denote by~$k^\star$ the color component with the
  largest~$|\lambda_\bfR^{\star,k}|$, one has:
  \begin{equation}
    \label{eqn:norm lambda max}
    \left\{
      \begin{array}{lclcl}
        \| \lambda \|_\infty & \egaldef &
        \displaystyle \max_{k=  R, G, B} \left| \lambda_\bfR^k \right| \leq
        \max_{k = R, G, B} \left| \lambda_\bfR^{\star, k} \right| =
        \left|\lambda_\bfR^{\star, k^\star} \right|,\\
        \| \lambda \|_\infty & = & \left| \lambda_\bfR^{\star, k^\star} \right|
        \quad \mbox{for the cutting } C = C_\bfR^{\star, k^\star}.
      \end{array}
    \right.
  \end{equation}
  This shows that~$C_\bfR^\star{\egaldef}C_\bfR^{\star,k^\star}$
  maximizes~$\|\lambda\|_\infty$ over the set~$T$ of all 2-partitions of the
  region~$\bfR$.
\end{description}

Since we have no theoretically funded way to choose a predefined set of
cuttings among which to select the refinement, the {\overallbest}
{\cuttingstrategy} will be the default option for image segmentation.

\subsubsection{Updating the partition $\calP$}

The algorithm computes the optimal cuttings~$C_\bfR^\star$ for each
region~$\bfR$ of the current partition~$\calP$ as indicated in the previous
section~\ref{sss:cutting_strategies}, and determines the region~$\bfR^\star$
with the largest {\em exact} refinement indicator.

The next partition~$\calP^\star$ is then obtained by splitting
region~$\bfR^\star$ according to cutting~$C_{\bfR^\star}^\star$, thus
increasing the number of regions by one.
The new segmented image~$c^\star$ is updated by formula~(\ref{eqn:cCi}) over
the region~$\bfR^\star$, and coincides with~$c$ over all other regions.

After~$n$ iterations, the {\em vector segmentation algorithm} adds~$n$ regions
to the initial partition.

\subsubsection{Optimal adaptive {\em vector segmentation} algorithm}
\label{sss:algo_vector_segmentation}

We summarize here the algorithm for the vector segmentation of an image~$d$.
Let~$I(C)$ denote the {\em best available} refinement indicators for the
ranking of cuttings~$C$ in the chosen family~$T$ (see
section~\ref{sss:cutting_strategies}).

\begin{description}
\item[Specify the {\cuttingstrategy}.]
  Choose one from:
  \begin{itemize}
  \item {\bestinafamily} \mbox{}
    (and then $I(C)={\Delta}J$ = exact indicators),
  \item {\overallbest} \mbox{}
    (and then $I(C)=\|\lambda\|_\infty$ = first order indicators).
  \end{itemize}

\item[Initialization.] \mbox{}
  \begin{enumerate}
  \item[-1.] \label{algo1:init}
    Define the initial partition~$\calP_1$ as the 1-region partition of the
    domain~$D$ of~$d$.
  \item[0.] \label{algo1:initmin}
    Compute the associated 1-region segmented image:
    \[
    c_1 = \calM_{\calP_1} \, \tilde{\calM}_{\calP_1} \, d
    \quad \Longleftrightarrow \quad
    \forall i \in I, (c_1)_i = \mbox{ mean value of } \{ d_j, j \in I \}.
    \]
  \end{enumerate}

\item[Iterations.]
  For $n{\geq}1$, do until the desired number of regions is attained:
  \begin{enumerate}
  \item[1.] \label{algo1:candidates}
    For all regions~$\bfR$ of the current partition~$\calP_n$,
    \begin{enumerate}
    \item[1a.] determine the best cutting~$C_\bfR^\star$ using the
      {\em best available} indicator~$I(C)$,
    \item[1b.] compute the {\em exact indicators}~${\Delta}J$ associated
      with~$C_\bfR^\star$ using~(\ref{eqn:exact_indicators_norm})
      and~(\ref{eqn:cCi}).
    \end{enumerate}
  \item[2.] Retain the region~$\bfR^\star$ whose best
    cutting~$C_{\bfR^\star}^\star$ has the largest {\em exact}
    indicator~${\Delta}J$.
  \item[3.] \label{algo1:keepbest}
    Define the next partition~$\calP_{n+1}$ by splitting~$\bfR^\star$
    according to cutting~$C_{\bfR^\star}^\star$.
  \item[4.] Define the next segmented
    image~$c_{n+1}=\calM_{\calP_{n+1}}\,\tilde{\calM}_{\calP_{n+1}}\,d$ by
    updating~$c_n$ on~$\bfR^\star$ according to the tentative segmented image
    on~$\bfR^\star$ obtained during step~1b.
  \end{enumerate}
\end{description}

\begin{remark}
When~$\calF$ is not the identity operator, one can use, for the
{\bestinafamily} {\cuttingstrategy}, the first order indicators in step~1a
instead of the exact ones if the cost of an exhaustive search with the exact
indicators is too high (see the adaptive parameterization
algorithm~\cite{ben:mri:08}).
\end{remark}

\begin{remark}
It is of course possible to start the algorithm from any partition~$\calP_1$
instead of a 1-partition.
\end{remark}

\subsection{Case of multiscalar segmentation}
\label{ss:case_multiscalar_segmentation}

In order to construct a multiscalar segmentation, one notice
that~(\ref{eqn:Jk_c}) is equivalent to
\begin{equation}
  \label{eqn:Jk_ck}
  \mbox{minimize } J^k (c^k) = \demi \sum_{i \in I} (d_i^k - c_i^k)^2
  \quad \mbox{for } k = R, G, B,
\end{equation}
under the same constraint that the number of regions of the
partitions~$\calP^k$ is ``small compared to the number of pixels of~$d$''.

This suggests to apply a scalar version of the previous optimal adaptive
segmentation algorithm {\em independently to each component} $k=R,G,B$.
If $\calP=(\calP^R,\calP^G,\calP^B)$ denotes the current multiscalar
partition, this will produce three {\em tentative optimal partitions}
$(\calP^{R\#},\calP^{G\#},\calP^{B\#})$.
It will then be necessary to choose one {\em \multiscalarstrategy} to decide
how to use this information to define the updated multiscalar segmentation
$\calP^\star=(\calP^{R\star}, \calP^{G\star},\calP^{B\star})$.

\subsubsection{Scalar optimal adaptive segmentation}
\label{sss:scalar_segmentation}

For a given color component~$k$, the {\em scalar algorithm} follows the steps
of section~\ref{sss:algo_vector_segmentation} with the following minor
adaptations.

Let~$\calP^k$ be the current partition for component~$k$.
The associated $k$-th component segmented image, solution of~(\ref{eqn:Jk_ck})
for the partition~$\calP^k$, is now given by
$c^k=\calM_{\calP^k}\tilde{\calM}_{\calP^k}d^k$, where
$\calM_{\calP^k}:(s_\bfR^k,\bfR\in\calP^k){\mapsto}c^k$ denotes the mapping
that attributes the segmented intensity~$s_\bfR^k$ to color component~$k$ on
each region~$\bfR$ of partition~$\calP^k$, and $\tilde{\calM}_{\calP^k}=
\left((\calM_{\calP^k})^T\calM_{\calP^k}\right)^{-1}(\calM_{\calP^k})^T$ its
least-squares pseudoinverse, which computes the mean value of the $k$-th
component over each region.

The refinement indicators associated with a cutting~$C$ which splits a
region~$\bfR$ of the current $k$-th component partition~$\calP^k$ into two
subregions~$\bfR_+$ and~$\bfR_-$ are (compare with
section~\ref{sss:ri_vector}):
\begin{description}
\item[Exact indicators:]
  \begin{equation}
    \label{eqn:exact_indicators_scalar}
    \Delta J^k = J^k (c^k) - J^k (c^k_C),
  \end{equation}
  where~$c^k$ is the $k$-th color component of the current segmented
  image~$c$, and~$c_C^k$ minimizes~$J^k$ defined in~(\ref{eqn:Jk_ck}) for the
  partition~$\calP_C^k$ obtained by splitting region~$\bfR$ of~$\calP^k$
  into~$\bfR_+$ and~$\bfR_-$.
  As in the vector case, $c_C^k$~is easily computed on~$\bfR$ by:
  \begin{equation}
    \label{eqn:ckCi}
    \left( c_C^k \right)_i = \left\{
      \begin{array}{lll}
        c_+^k & \mbox{for } x_i \in \bfR_+ &
        \mbox{where~$c_+^k$ = mean value of } \{ d_j^k, x_j \in \bfR_+ \},\\
        c_-^k & \mbox{for } x_i \in \bfR_- &
        \mbox{where~$c_-^k$ = mean value of } \{ d^k_j, x_j \in \bfR_- \}.
      \end{array}
    \right.
  \end{equation}
  
\item[First order indicators:]
  \begin{equation}
    \label{eqn:eqn:first_order_indicators_scalar}
    | \lambda^k | \mbox{ with } \lambda^k
    \mbox{ given by (\ref{eqn:first_order_indicators_RGB})}.
  \end{equation}
  Of course, exact and first order indicators are still linked by
  relation~(\ref{eqn:identite_remarquable}).
\end{description}

In the absence of specific information for each component, we have chosen to
use {\em the same set}~$T$ of tentative cuttings for each color component~$k$.
So the {\cuttingstrategies} are (compare to
section~\ref{sss:cutting_strategies}):
\begin{description}
\item[{\Bestinafamily}.]
  $T$~is a {\em predefined family} of cuttings, for example the set of vertical
  and horizontal cuttings.
  With this strategy, the best cutting~$C_\bfR^\star$ of a region~$\bfR$
  of~$\calP^k$ for the {\em exact indicator}~$\Delta J^k$ needs to be
  computed by an exhaustive search.
  
\item[{\Overallbest}.]
  $T$~is made of {\em all possible 2-partitions} of the region~$\bfR$.
  Here the best cutting~$C_\bfR^\star$ for the
  {\em first order indicators}$|\lambda^k|$ can be determined analytically by
  formula~(\ref{eqn:first_order_indicators_opt_RGB}).
\end{description}

So after one iteration of the scalar algorithm has been performed on each
component~$k$, the following quantities are available:
\begin{eqnarray}
  \label{eqn:pkdieze}
  \calP^{k \#} & = &
  \mbox{tentative optimal partition for color component~$k$}, \\
  \label{eqn:djkdieze}
  \Delta J^{k \#} & = &
  \mbox{exact decrease of the minimum of~$J^k$ when~$\calP^k$ is replaced
    by~$\calP^{k \#}$}.
\end{eqnarray}

\subsubsection{Updating the partition
  $\calP=(\calP^R,\calP^G,\calP^B)$: {\multiscalarstrategy}}
\label{sss:ms_strategy}

We describe now three {\em \multiscalarstrategies} which can be applied to
update the current multiscalar partition~$\calP$ once the
results~(\ref{eqn:pkdieze}) and~(\ref{eqn:djkdieze}) of the scalar algorithm
have been obtained for each color component:
\begin{description}
\item[{\Bestcomponentonly}.]
  This strategy determines the color component~$k^\star$ whose tentative
  refinement produces the largest decrease~$\Delta J^{k\#}$, and applies the
  corresponding optimal partition~$\calP^{k^\star\#}$ to component~$k^\star$
  only.
  When~$k^\star=R$ (red) for example, the current and updated multiscalar
  partitions are:
  \begin{equation}
    \label{eqn:best_component_only}
    \calP = (\calP^R, \calP^G, \calP^B),
    \quad \calP^\star = (\calP^{R \#}, \calP^G, \calP^B).
  \end{equation}
  After $n$~iterations, the {\bestcomponentonly} {\multiscalarstrategy}
  adds $n$~regions to the initial multiscalar partition.
  Of course, the vector partition corresponding to the resulting multiscalar
  partition, obtained by superimposing the three scalar partitions, may contain
  many more regions.
  
\item[{\Bestcomponentforeach}.]
  In this strategy, the algorithm refines each component~$k$ according to the
  corresponding tentative optimal partition~$\calP^{k\#}$ determined by the
  scalar algorithms.
  The current and updated multiscalar partitions are then:
  \begin{equation}
    \label{eqn:best_component_for_each}
    \calP = (\calP^R, \calP^G, \calP^B),
    \quad \calP^\star = (\calP^{R \#}, \calP^{G \#}, \calP^{B \#}).
  \end{equation}
  After $n$~iterations, the {\bestcomponentforeach} {\multiscalarstrategy} has
  added $3n$~regions to the initial partition ($n$~regions per component).
  Again, the vector partition corresponding to the resulting multiscalar
  partition, obtained by superimposing the three scalar partitions, may contain
  many more regions.
  
  The {\bestcomponentforeach} {\multiscalarstrategy}
  amounts to apply the scalar refinement indicator algorithm to each component
  separately, i.e. to use the {\em tensor product} of the three algorithms.
  
\item[{\Combinebestcomponents}.]
  Here the algorithm selects the three optimal
  cuttings~$C_{\bfR^\star}^{k\star}$ corresponding to the best tentative
  segmentation for each component $k=R,G,B$, and applies {\em all of them} to
  {\em all components}.
  This makes sense only if the three current component partitions coincide, and
  produces a partition~$\calP^{RGB\#}$, which is nothing but the
  superimposition of the three tentative partitions~$\calP^{R\#}$,
  $\calP^{G\#}$ and~$\calP^{B\#}$.
  This partition is used to update all three components:
  \begin{equation}
    \label{eqn:combine_best_components}
    \calP = (\calP^R, \calP^G, \calP^B),
    \quad \calP^\star = (\calP^{RGB \#}, \calP^{RGB \#}, \calP^{RGB \#}).
  \end{equation}
  So the {\combinebestcomponents} {\multiscalarstrategy} constructs in fact a
  sequence of {\em vector segmentations} of the original image~$d$.
  
  In contrast with the vector segmentation algorithm of
  section~\ref{sss:algo_vector_segmentation}, the {\combinebestcomponents}
  {\multiscalarstrategy} adds at each iteration between~$3$ and~$2^{3}-1$
  regions to the current vector partition.
  So the number of regions may increase rapidly, which is not compatible with
  the idea of adaptive parameterization.
\end{description}

Note that all {\multiscalarstrategies} are equivalent in the scalar case.

\subsubsection{The multiscalar adaptive segmentation algorithm}
\label{sss:algo_multiscalar_segmentation}

We summarize here the adaptation of the algorithm of
section~\ref{sss:algo_vector_segmentation} to the determination of
multiscalar segmentations.

\begin{description}
\item[Specify the {\cuttingstrategy}.]
  Choose one from:
  \begin{itemize}
  \item {\bestinafamily} \mbox{}
    (best available = exact indicators $I^k(C)=\Delta J^k$),
  \item {\overallbest} \mbox{}
    (best available = first order indicators $I^k(C)=|\lambda^k|$).
  \end{itemize}
  
\item[Specify the {\multiscalarstrategy}.]
  Choose one from:
  \begin{itemize}
  \item {\bestcomponentonly},
  \item {\bestcomponentforeach},
  \item {\combinebestcomponents}.
  \end{itemize}

\item[Initialization.]
  Define $\calP_1=(\calP_1^R,\calP_1^G,\calP_1^B)$ and
  $c_1=(c_1^R,c_1^G,c_1^B)$ by:
  \begin{description}
  \item[Scalar initializations.]
    For each component $k=R,G,B$:
    \begin{enumerate}
    \item[-1.] \label{algo2:init}
      Define the~$\calP_1^k$ as the 1-region partition of the domain~$D$
      of~$d$.
    \item[0.] \label{algo2:initmin}
      Compute the associated 1-region segmented {\em component image}:
      \[
      c_1^k = \calM_{\calP_1^k} \, \tilde{\calM}_{\calP_1^k} \, d^k
      \, \Longleftrightarrow \,
      \forall i \in I, (c_1^k)_i = \mbox{ mean value of } \{ d_j^k, j \in I \}.
      \]
    \end{enumerate}
  \end{description}
  
\item[Iterations]
  For $n{\geq}1$, do until the desired number of regions is attained:
  \begin{description}
  \item[Scalar iterations.]
    For each component $k=R,G,B$, do (section~\ref{sss:scalar_segmentation}):
    \begin{enumerate}
    \item[1.] \label{algo2:candidates}
      For all regions~$\bfR$ of the current partition~$\calP_n^k$,
      \begin{enumerate}
      \item[1a.] determine the best cutting~$C_\bfR^\star$ using the
        {\em best available} indicator~$I(C)$,
      \item[1b.] compute the {\em exact indicators}~${\Delta J}^k$ associated
        with~$C_\bfR^\star$ using~(\ref{eqn:exact_indicators_scalar})
        and~(\ref{eqn:ckCi}).
      \end{enumerate}
    \item[2.] \label{algo2:keepbest}
      Retain the region~$\bfR^\star$ whose best
      cutting~$C^{k\star}_{\bfR^\star}$ has the largest {\em exact}
      indicator~${\Delta J}^k$.
    \item[3.] Define the next tentative partition~$\calP^{k\#}$ of color
      component~$k$ by splitting the region~$\bfR^\star$ of~$\calP^k$
      according to cutting~$C^{k\star}_{\bfR^\star}$.
    \end{enumerate}
    
  \item[Update the multiscalar segmentation.] \mbox{}
    \begin{enumerate}
    \item[4.] Compute the next multiscalar partition~$\calP_{n+1}$ according to
      the chosen {\multiscalarstrategy} (section~\ref{sss:ms_strategy}).
    \item[5.] Define the next multiscalar segmented image~$c_{n+1}$:
      \begin{itemize}
      \item {\bestcomponentonly}:
        let~$k^\star$ be the component retained at step~4.
        It suffices then to update~$c_n^{k^\star}$ on the region~$\bfR^\star$
        determined at step~2 according to the tentative segmented
        image~(\ref{eqn:ckCi}) determined at step~1b.
      \item {\bestcomponentforeach}:
        each component~$k$ of~$c_n^k$ needs to be updated on the
        corresponding region~$\bfR^\star$ determined at step~2 according to
        the tentative segmented image~(\ref{eqn:ckCi}) obtained during
        step~1b.
      \item {\combinebestcomponents}:
        here~$c_{n+1}$ is obtained by updating~$c_n$ on the~1 to~3 regions
        $\bfR\in\calP_n$ which have been split to define~$\calP_{n+1}$.
      \end{itemize}
    \end{enumerate}
  \end{description}
\end{description}

\begin{remark}
It is of course possible to start the algorithm from any multiscalar
partition~$\calP_1$ instead of three 1-partitions, with the restriction
that~$\calP_1$ has to be a {\em vector partition} in case the
{\combinebestcomponents} {\multiscalarstrategy} is used.
\end{remark}

\begin{remark}
Note that for both the {\em vector} and {\em multiscalar} segmentations, the
expressions to compute are simple closed-form formulas that are local to
regions.
Furthermore, in steps~1a and~1b the best cutting for each region and its
associated exact indicator are invariant as long as the region is not split.
Hence, an iteration of the algorithm only amounts to compute step~1 for the
two new subregions (of the retained region~$\bfR^\star$) that appeared in the
partition~$\calP_n$ at the previous iteration.
Taking into account this optimization, it is obvious to derive a linear version
of the algorithm (in terms of the number of iterations), which is a salient
aspect of this approach.
\end{remark}


\section{Numerical results}
\label{s:numerical_results}

We illustrate in this section the numerical behavior of algorithms of
sections~\ref{sss:algo_vector_segmentation}
and~\ref{sss:algo_multiscalar_segmentation} for the adaptive vector and
multiscalar segmentations.
The image~$d$ to be segmented is the picture of a colorful
sculpture~\cite{grosbec}.
Its definition is $400\times533$ (i.e., 213,200 pixels).

Note that the original image~$d$ is itself trivially a segmented image, with a
number of regions equal, for a vector segmentation, to the number of distinct
colors---here 55,505---, or equal, for each component for a multiscalar
segmentation, to the number of distinct levels in the component---here the
maximum of 256 levels for a 24-bit color image is reached for all components.
Hence, if no constraint is set on the number of regions, all above adaptive
segmentation algorithms will converge towards the original image~$d$.
But, though regions are introduced in our algorithms with parsimony at each
iteration, there is for the time no mathematical proof that they will produce a
segmented image~$c$ identical to~$d$ as soon as the number of regions in $c$
becomes equal to that of the trivial segmentation of $d$.
But~$c$ will for sure end up being identical to $d$ at some point before the
number of its regions equals the number of its pixels.

This convergence property is important as it ensures that by properly stopping
the algorithm, one can control the compromise between the number of regions of
the segmentation and the resemblance of the segmented image to the original
one.

The software used to compute these numerical examples is still under
development, but it will be available soon at
{\tt http://refinement.inria.fr/}.

\subsection{Vector segmentation}
\label{ss:vector_segmentation_numerical}

We test first the vector segmentation
algorithm~\ref{sss:algo_vector_segmentation} with two {\cuttingstrategies}.
Here the algorithm adds exactly one region at each iteration, so it is easy to
specify a stopping criterion providing a segmentation into exactly~$n$ regions
with uniform colors.

\subsubsection{{\Bestinafamily} {\cuttingstrategy}}

\begin{figure}[htb]
  \begin{center}
    (a)
    \includegraphics[width=\quartps\textwidth]{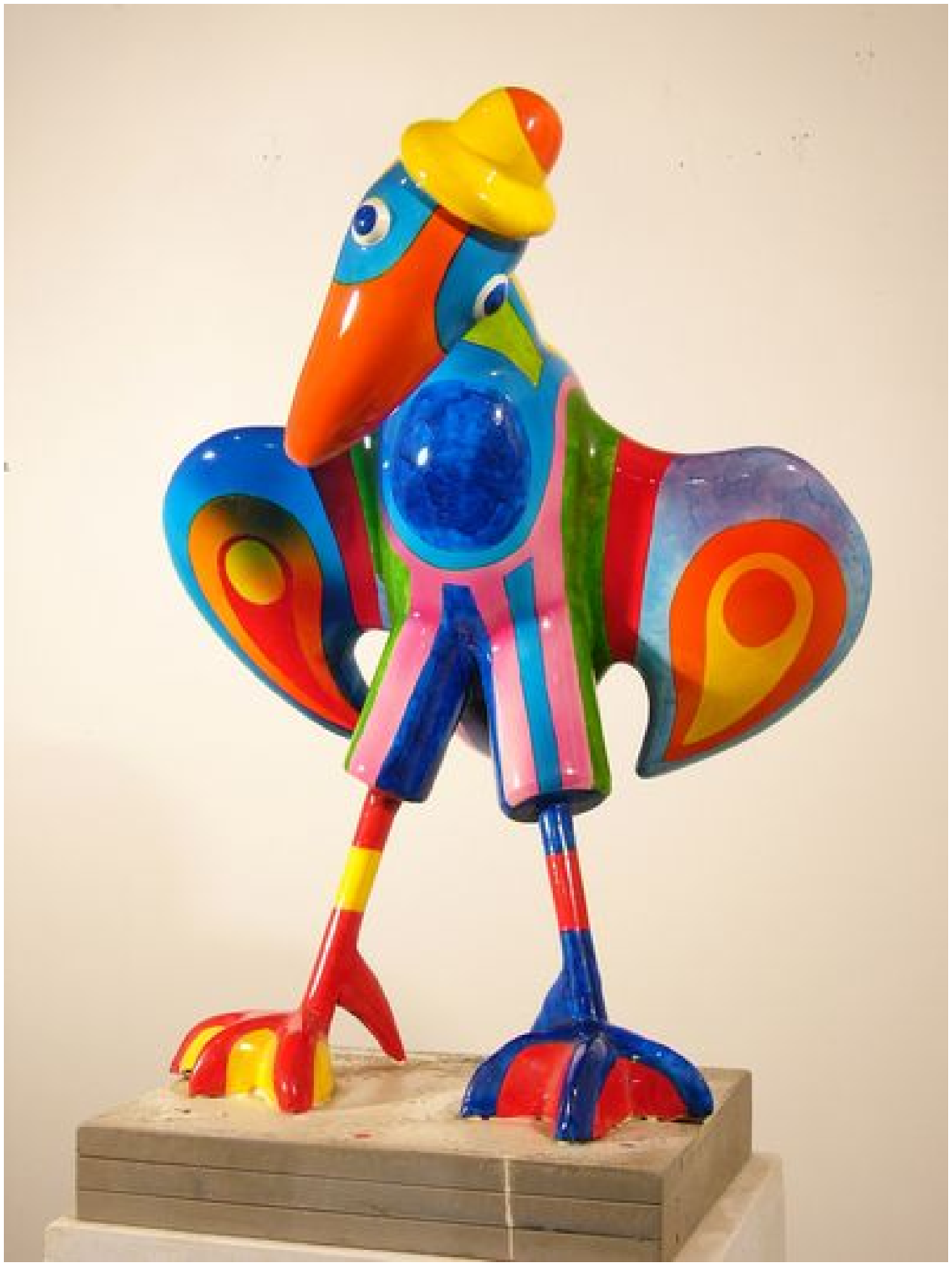}
    (b)
    \includegraphics[width=\quartps\textwidth]{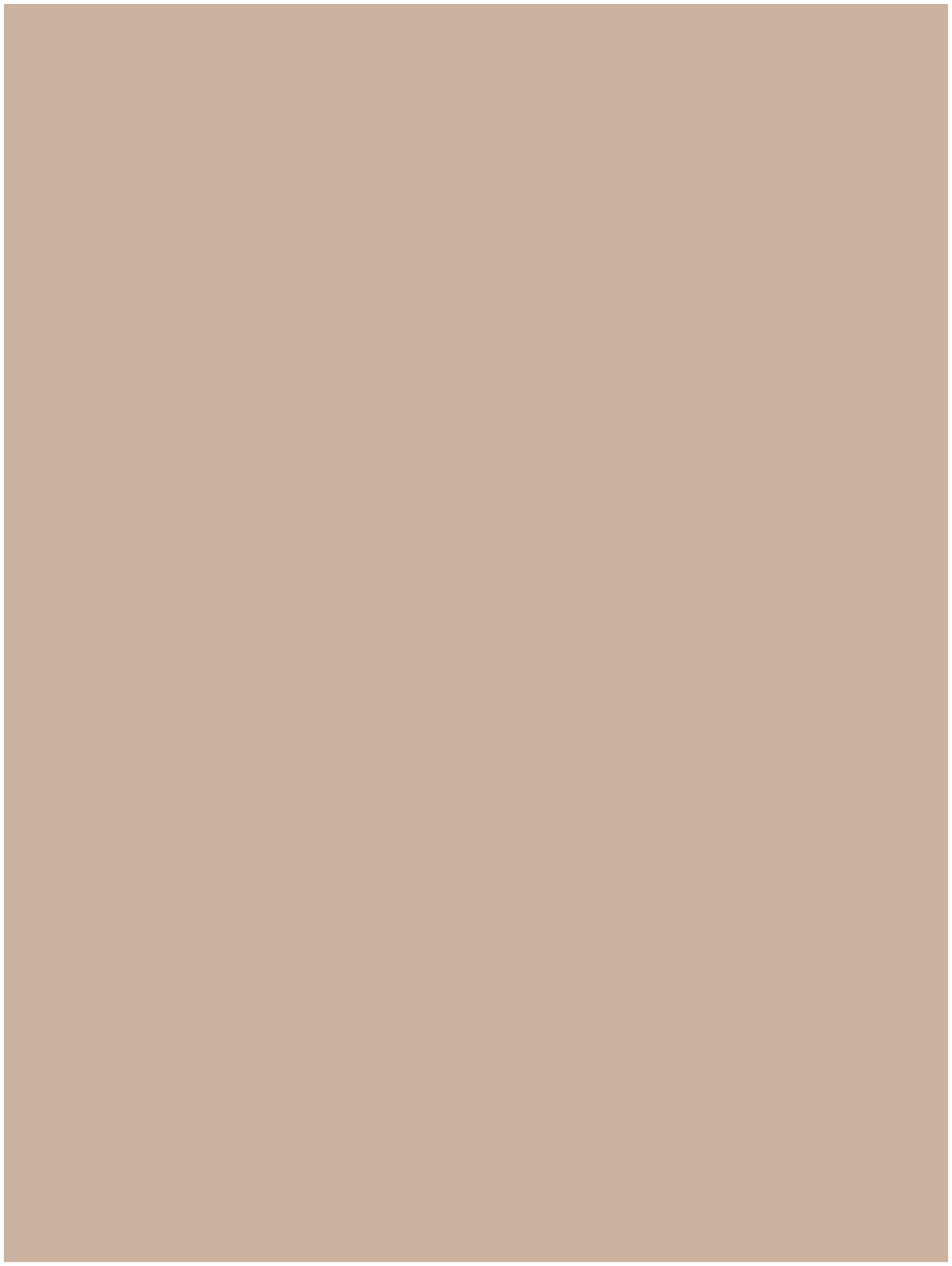}
    (c)
    \includegraphics[width=\quartps\textwidth]{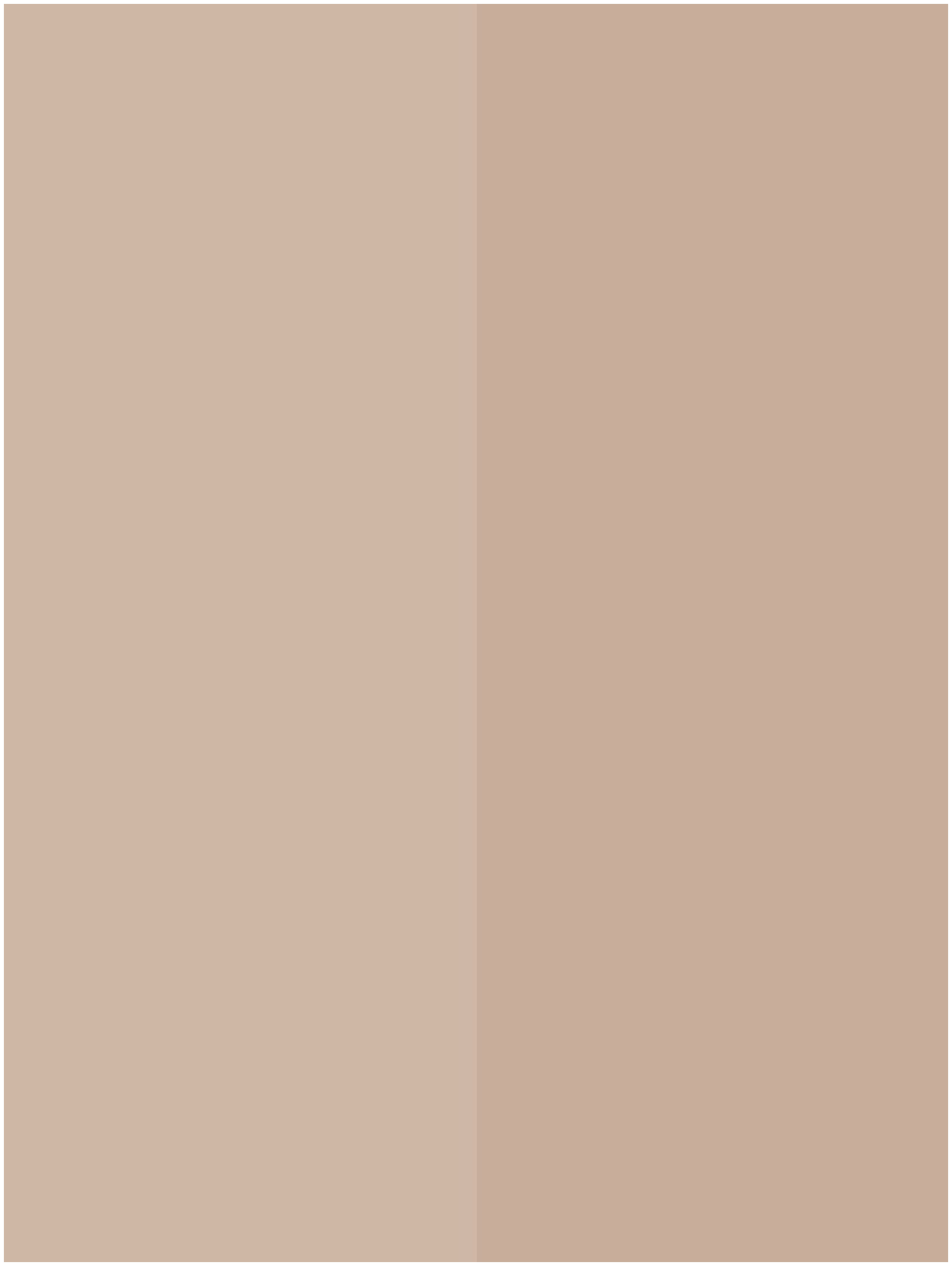} \\
    (d)
    \includegraphics[width=\quartps\textwidth]{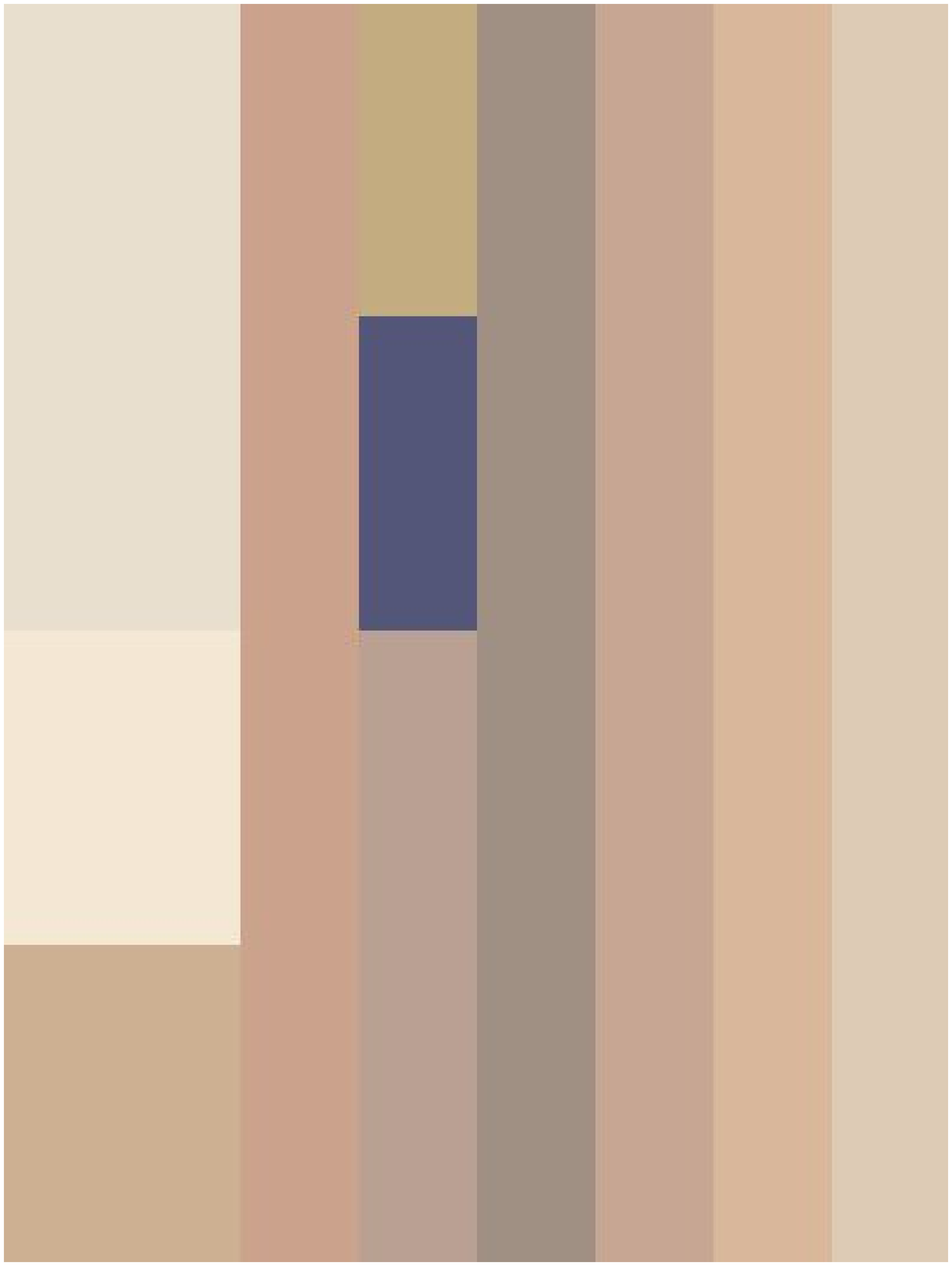}
    (e)
    \includegraphics[width=\quartps\textwidth]{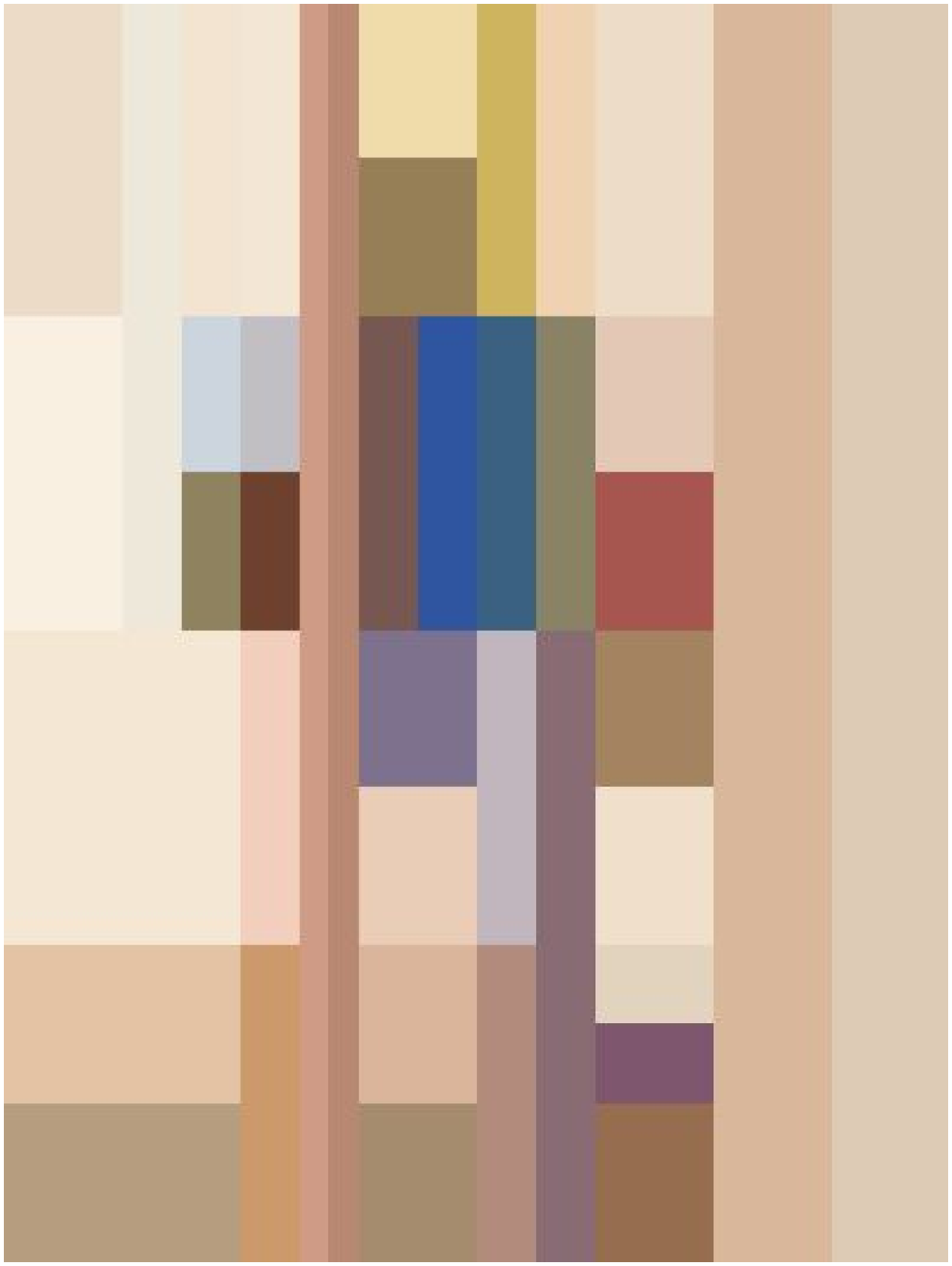}
    (f)
    \includegraphics[width=\quartps\textwidth]{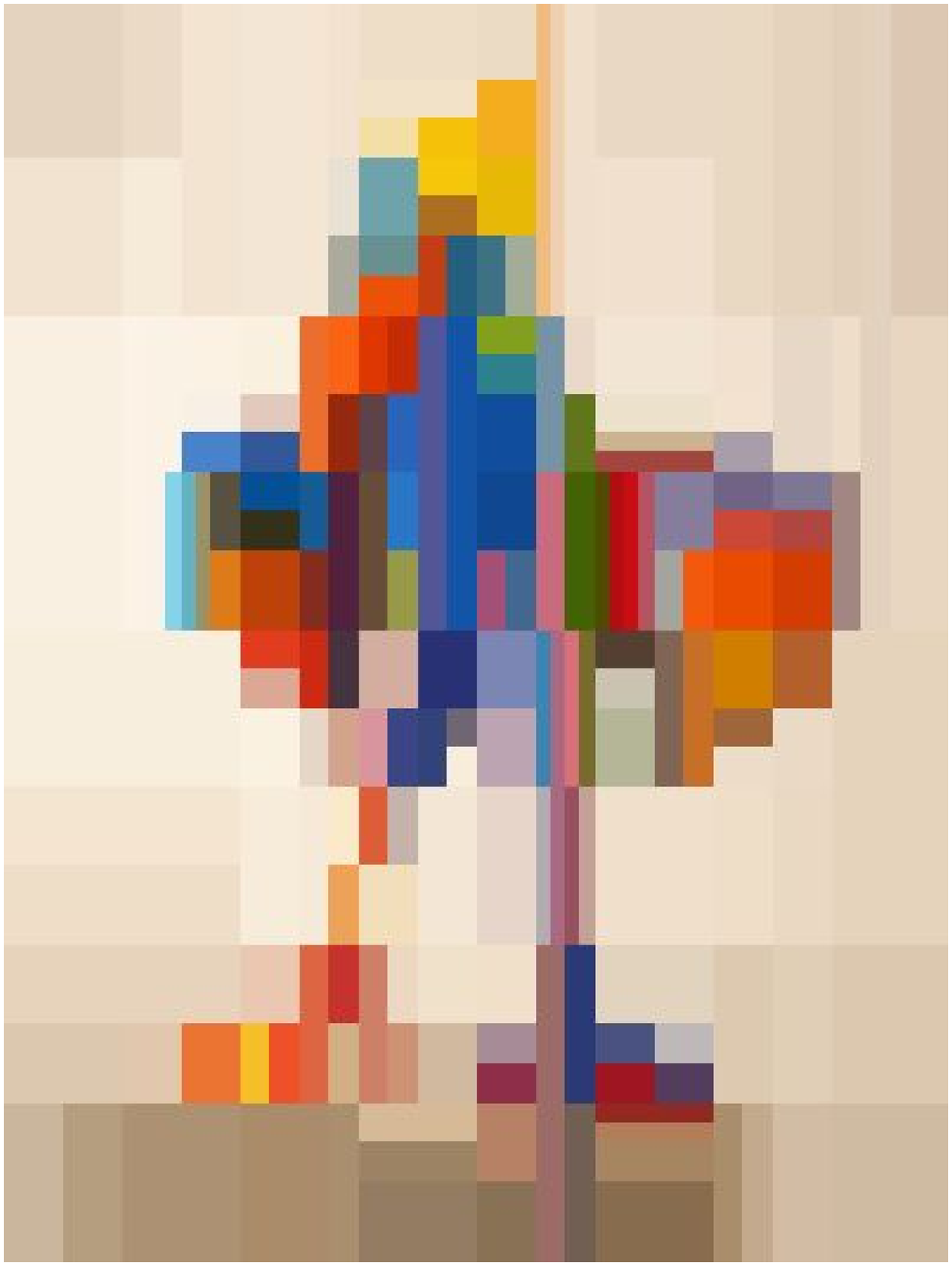}
    \caption{{\em Vector segmentation} with the {\bestinafamily}
      {\cuttingstrategy}.
      (a) is the image data to segment.
      Then, for each displayed image, we give the iteration number~$n$,
      the number of regions~$n_{\rm vr}$, and the percentage of explained
      data~$\tau$.\protect\\
      (b) $n=0$, $n_{\rm vr}=1$, $\tau=65.0\%$;\protect\\
      (c) $n=1$, $n_{\rm vr}=2$, $\tau=65.1\%$;\protect\\
      (d) $n=10$, $n_{\rm vr}=11$, $\tau=68.7\%$;\protect\\
      (e) $n=40$, $n_{\rm vr}=41$, $\tau=73.4\%$;\protect\\
      (f) $n=200$, $n_{\rm vr}=201$, $\tau=82.2\%$.}
    \label{fig:bestinafamily}
  \end{center}
\end{figure}

The images shown in Figure~\ref{fig:bestinafamily} are obtained with a quite
scarce predefined family of vertical and horizontal cuttings: regions are
rectangles, and the only cuttings to be tested in each region are its
horizontal and vertical divisions into two equal parts.
For instance, the first cutting is vertical as shown
in~Figure~\ref{fig:bestinafamily}(c), and all regions
in~Figures~\ref{fig:bestinafamily}(d), \ref{fig:bestinafamily}(e)
and~\ref{fig:bestinafamily}(f) are rectangular.
The exact indicators are used to rank all cuts.

The adaptive character of the algorithm is clear on images
Figures~\ref{fig:bestinafamily}(d), \ref{fig:bestinafamily}(e)
and~\ref{fig:bestinafamily}(f) where regions are extremely irregular.
At each iteration, the cutting providing the largest decrease of the cost
function is selected.

The convergence of the algorithm is very slow: the image
Figure~\ref{fig:bestinafamily}(f) obtained after 200~iterations is made of
200~monochrome rectangles, but it bears only a rough resemblance to the data.
So it seems that the {\bestinafamily} {\cuttingstrategy} is not well
suited for segmentation purpose---even though it remains a very good strategy
for the general problem of vector parameter estimation.

\subsubsection{{\Overallbest} {\cuttingstrategy}}

\begin{figure}[htb]
  \begin{center}
    (a)
    \includegraphics[width=\quartps\textwidth]{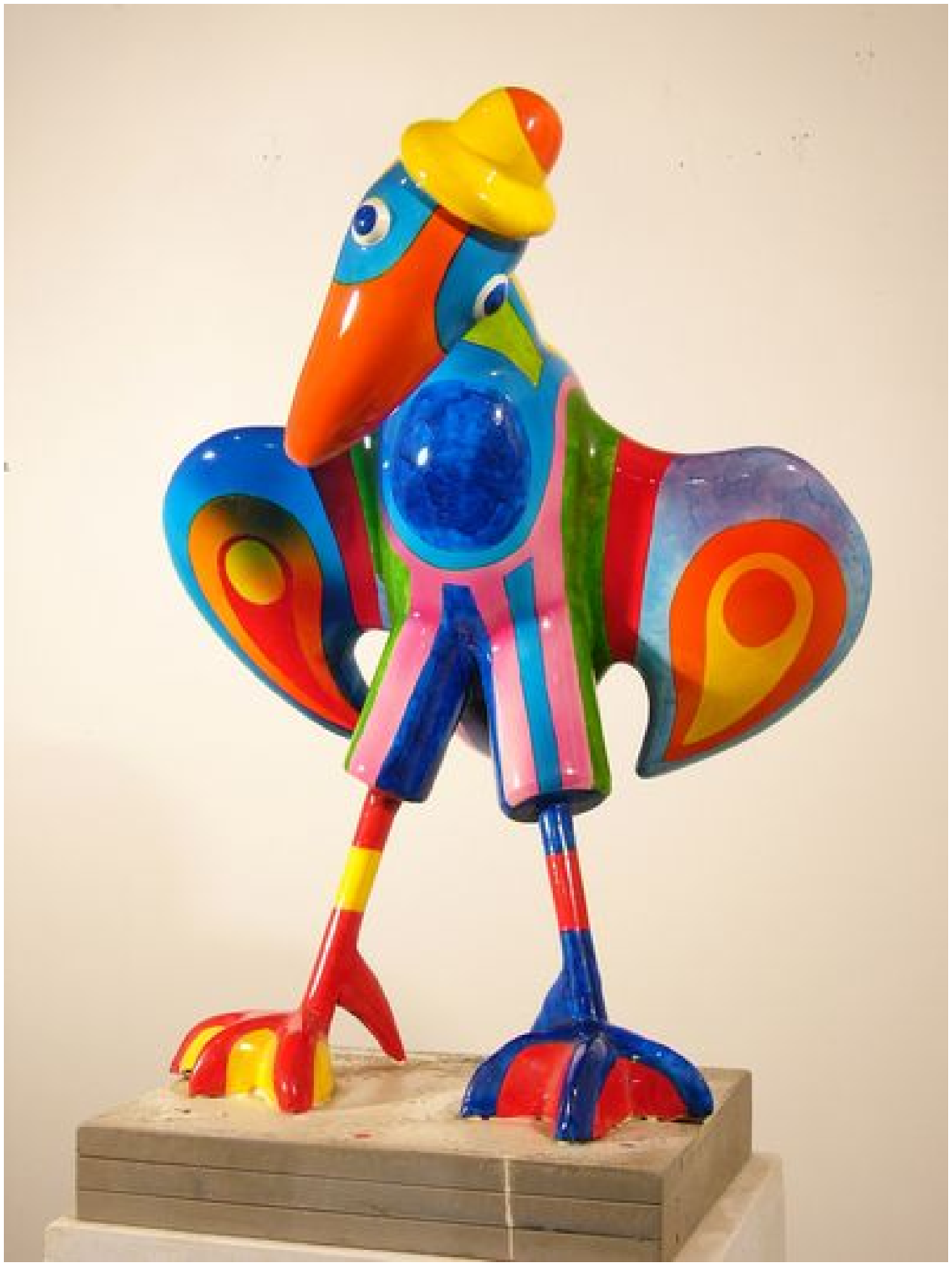}
    (b)
    \includegraphics[width=\quartps\textwidth]{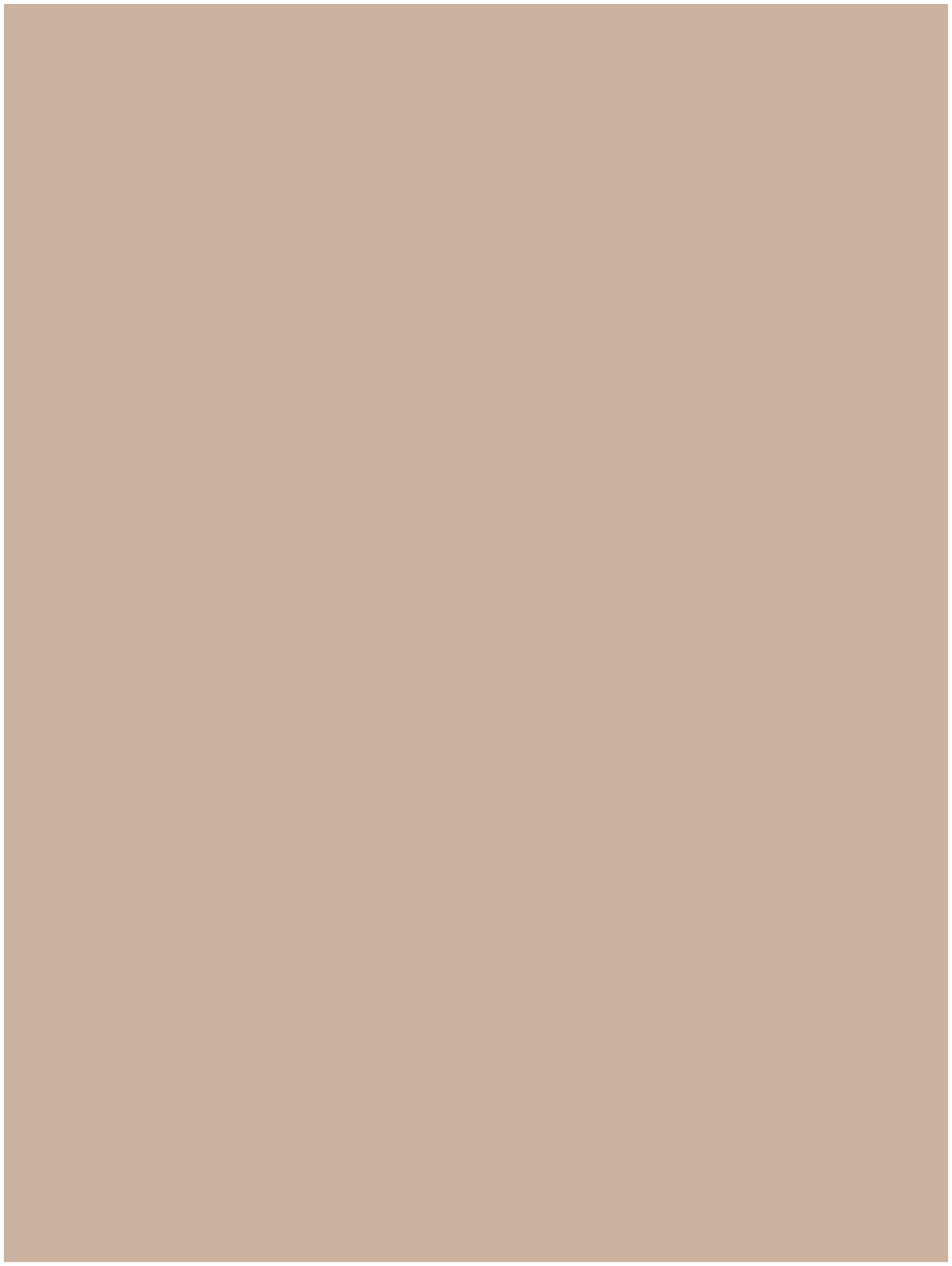}
    (c)
    \includegraphics[width=\quartps\textwidth]{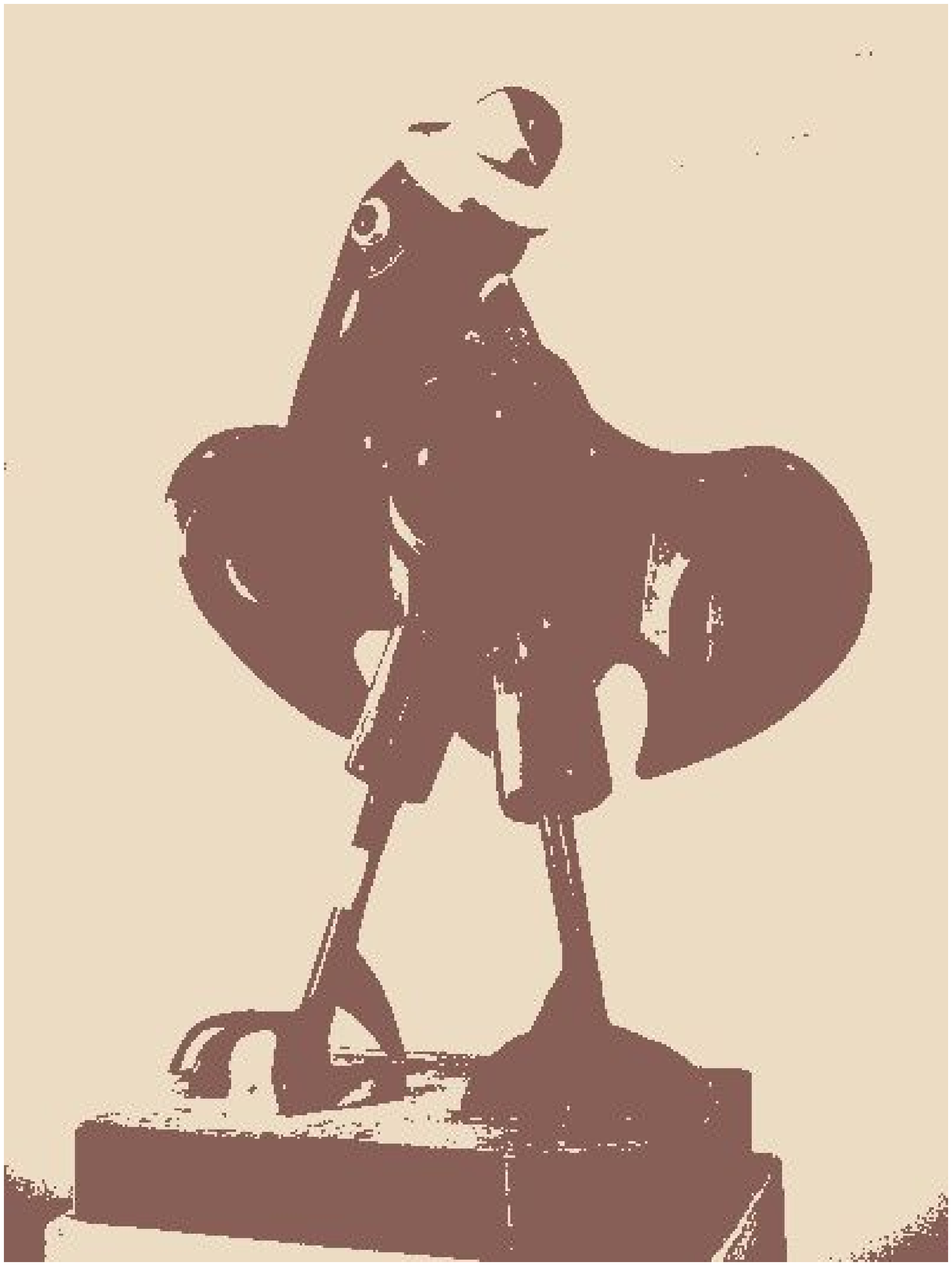}
    (d)
    \includegraphics[width=\quartps\textwidth]{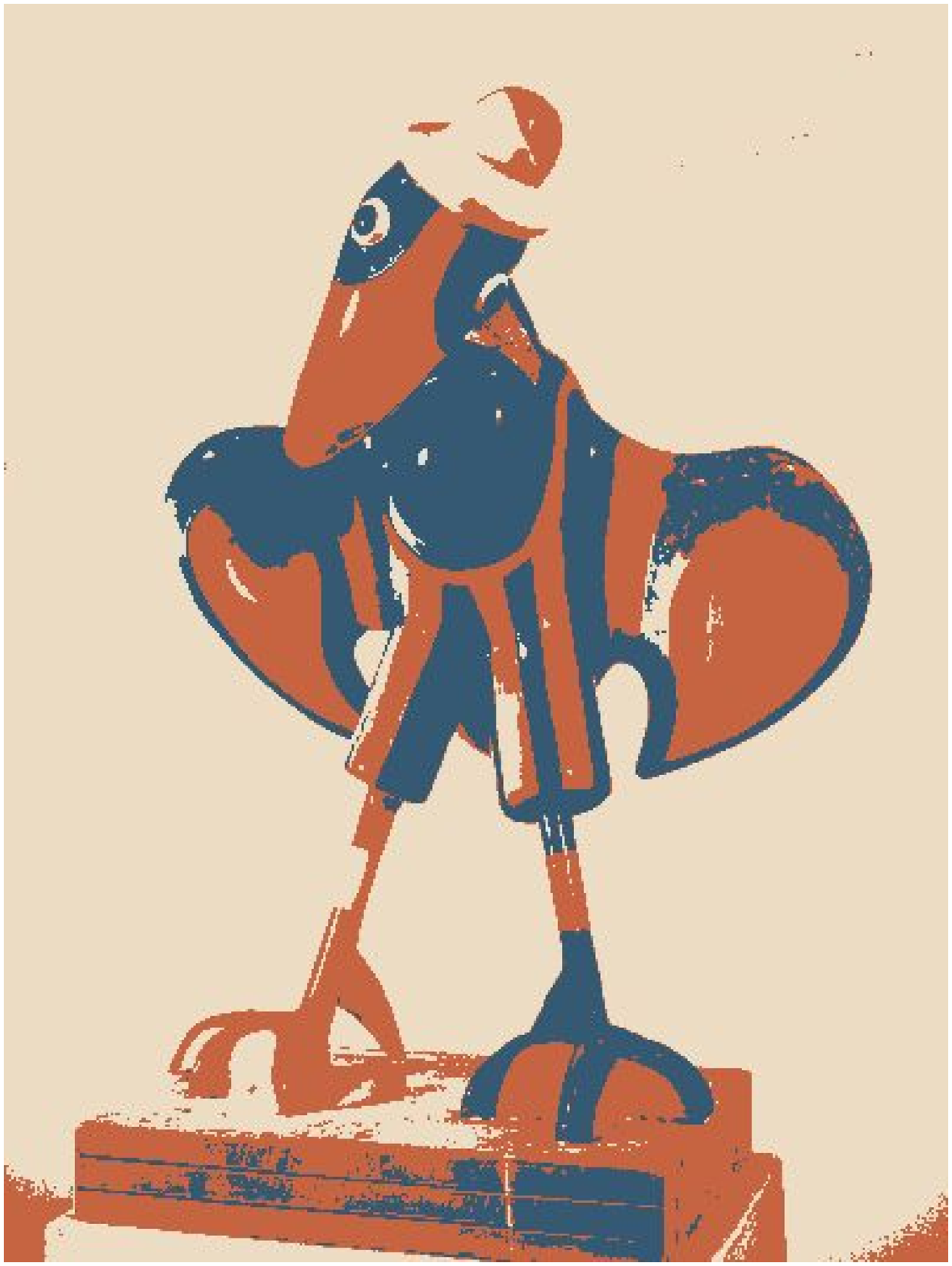} \\
    (e)
    \includegraphics[width=\quartps\textwidth]{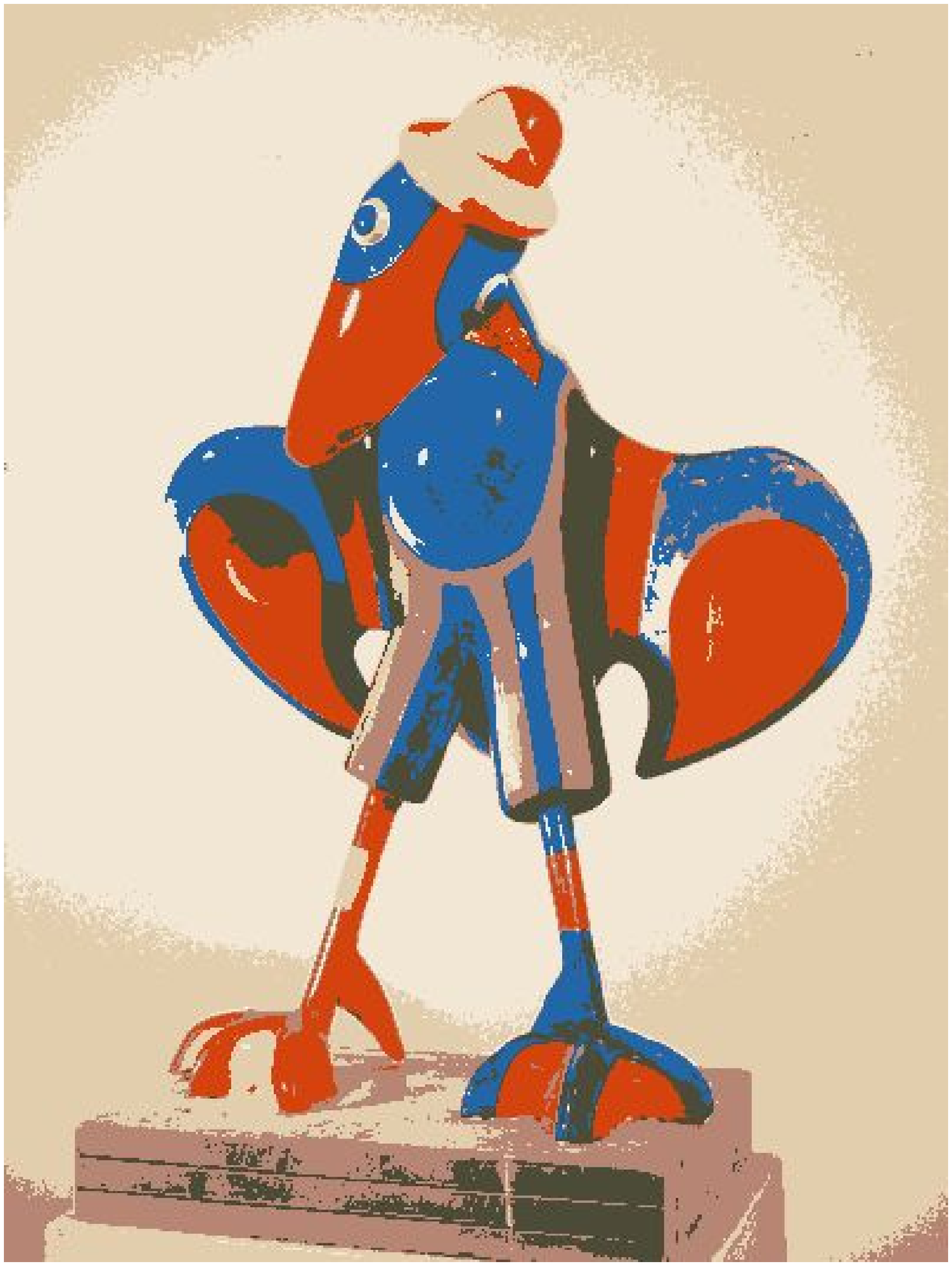}
    (f)
    \includegraphics[width=\quartps\textwidth]{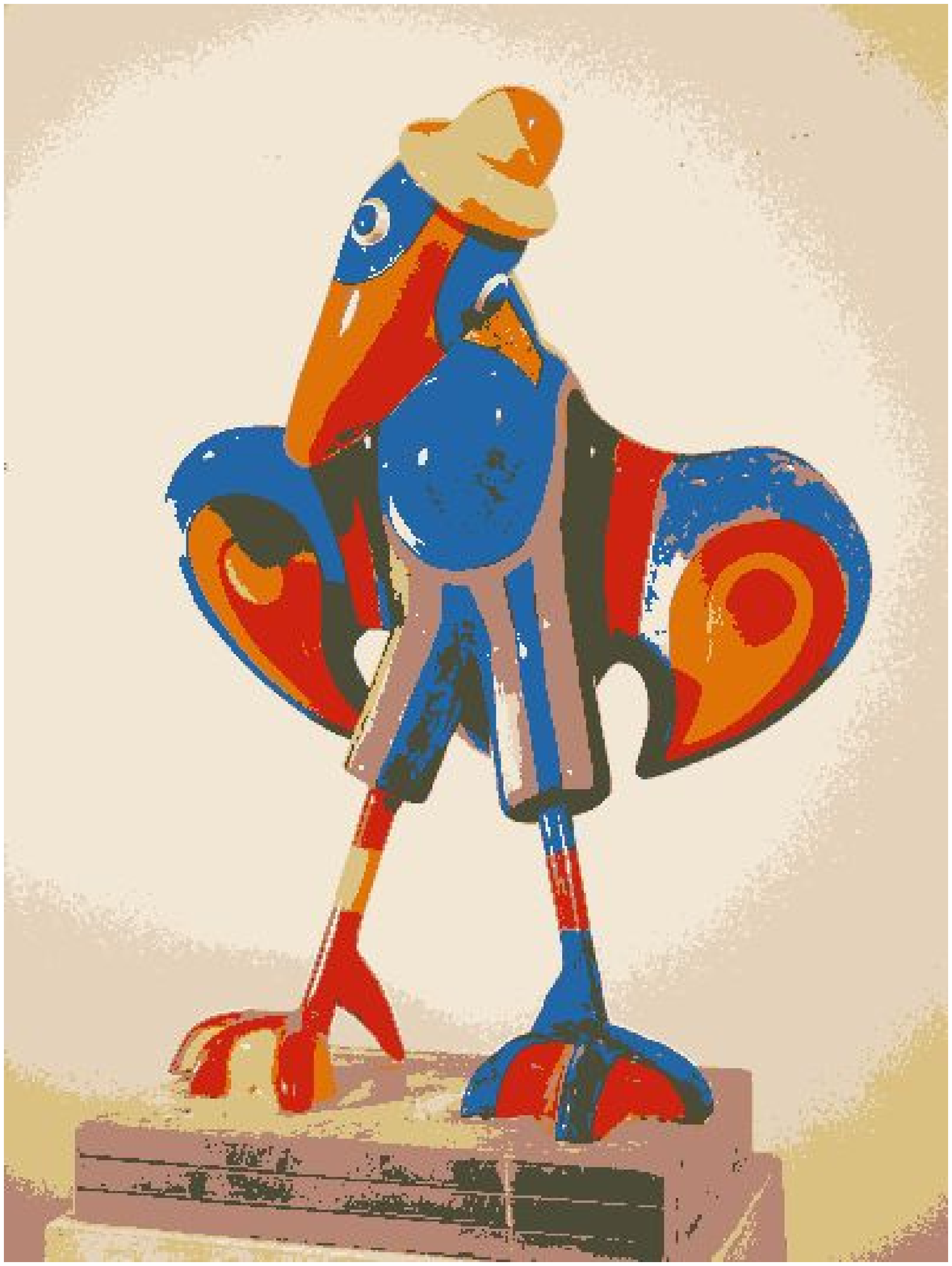}
    (g)
    \includegraphics[width=\quartps\textwidth]{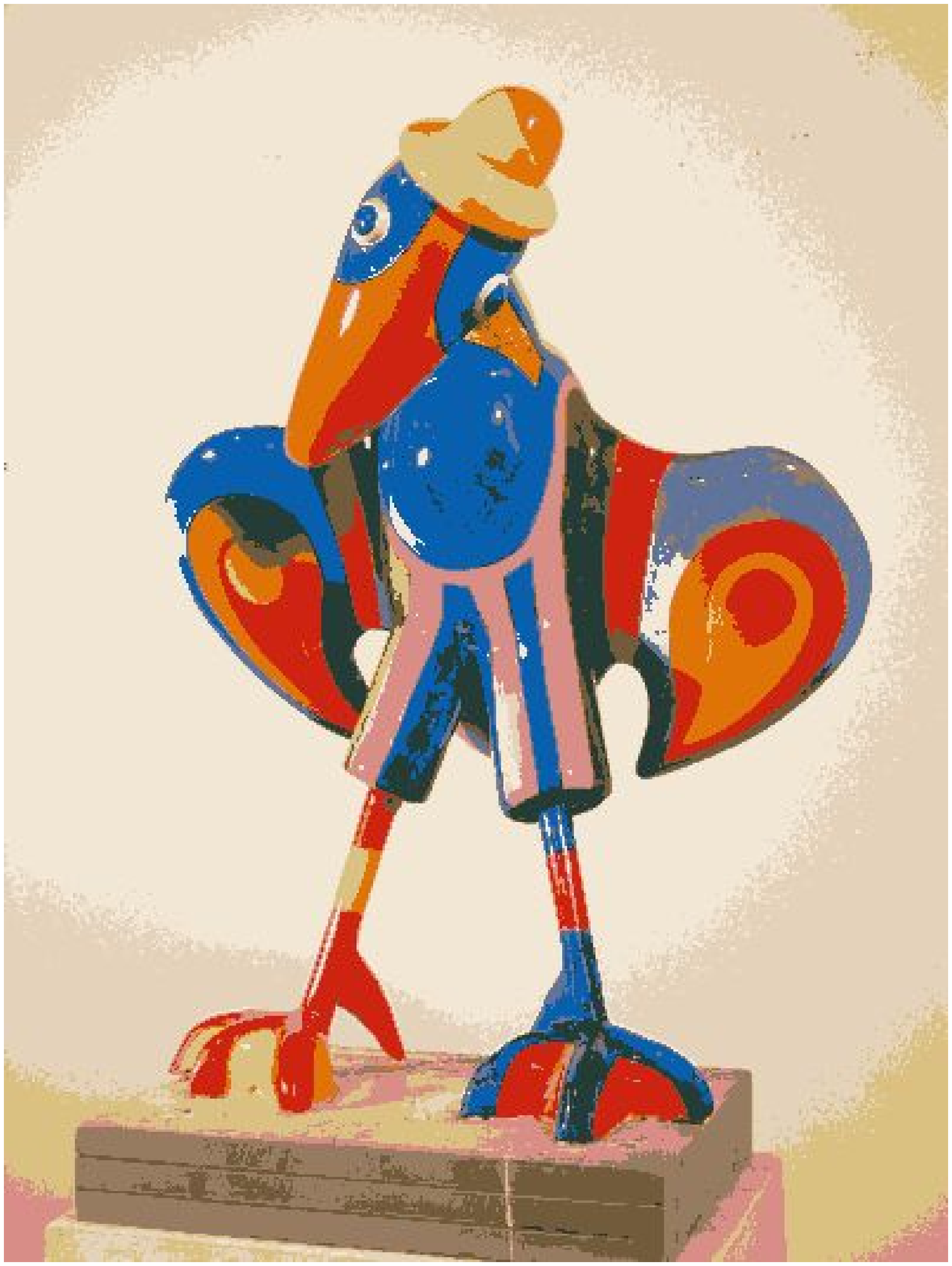}
    (h)
    \includegraphics[width=\quartps\textwidth]{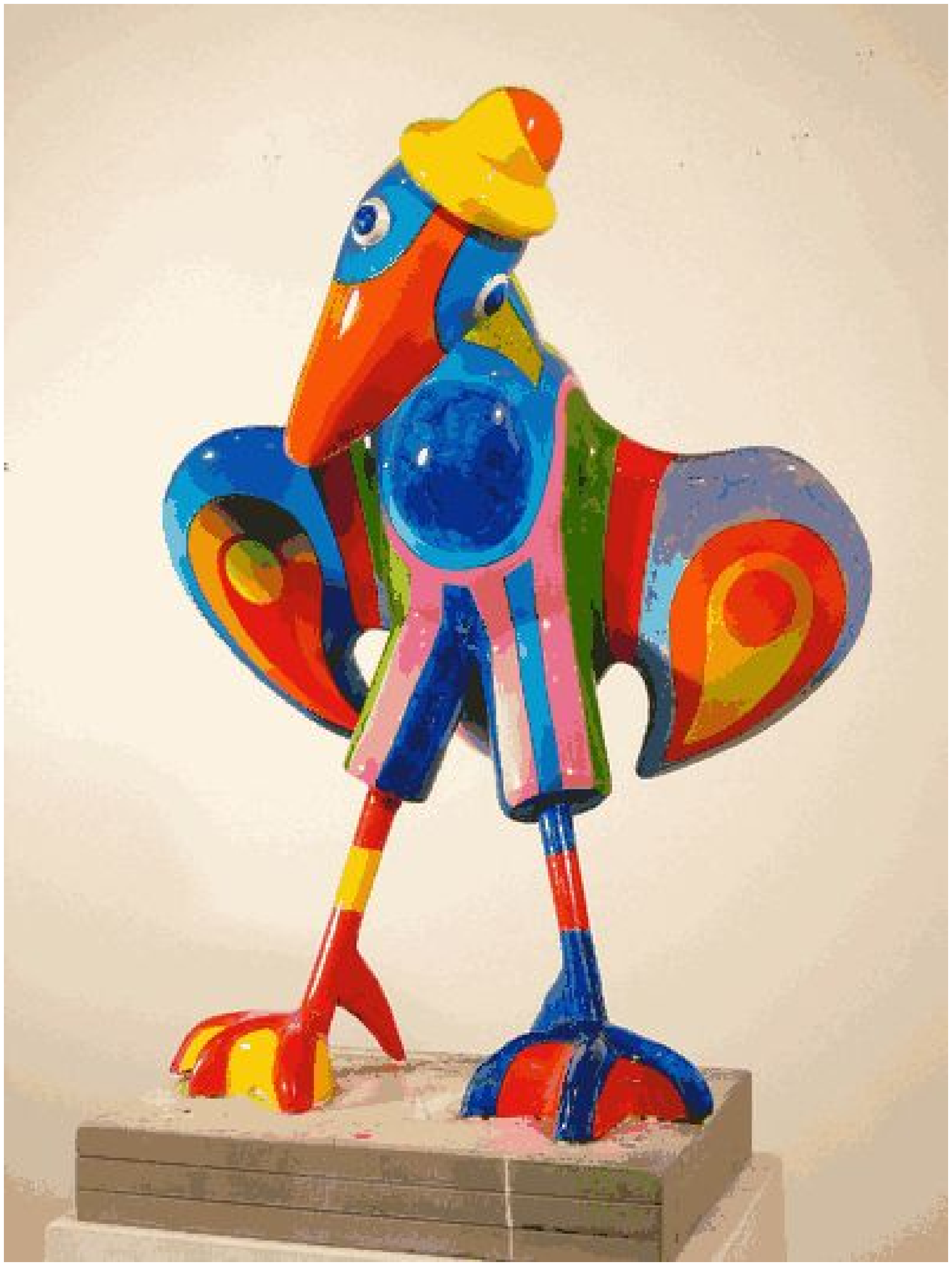}
    \caption{{\em Vector segmentation} with the {\overallbest}
      {\cuttingstrategy}.
      (a) is the image data to segment.
      Then, for each displayed image, we give the iteration number~$n$,
      the number of regions~$n_{\rm vr}$, and the percentage of explained
      data~$\tau$.\protect\\
      (b) $n=0$, $n_{\rm vr}=1$, $\tau=65.0\%$;\protect\\
      (c) $n=1$, $n_{\rm vr}=2$, $\tau=77.9\%$;\protect\\
      (d) $n=2$, $n_{\rm vr}=3$, $\tau=82.1\%$;\protect\\
      (e) $n=5$, $n_{\rm vr}=6$, $\tau=86.8\%$;\protect\\
      (f) $n=7$, $n_{\rm vr}=8$, $\tau=88.1\%$;\protect\\
      (g) $n=10$, $n_{\rm vr}=11$, $\tau=89.6\%$;\protect\\
      (h) $n=40$, $n_{\rm vr}=41$, $\tau=94.8\%$.}
    \label{fig:overallbest}
  \end{center}
\end{figure}

We present in Figure~\ref{fig:overallbest} the vector segmentation
of the same image with the {\overallbest} {\cuttingstrategy}.
Here the shape of the regions are perfectly adapted to the image data: the
image is already recognizable with only three regions,
see Figure~\ref{fig:overallbest}(d), and the algorithm is able to
explain 86.8\%~of the data with as few as six color regions, see
Figure~\ref{fig:overallbest}(e).
But it takes 41~color regions to explain up to 94.8\%~of the data.
In any case, the intermediate images generated by the algorithm offer a large
choice of segmented images to choose from.

It is this configuration of the algorithm (vector segmentation with
{\overallbest} {\cuttingstrategy}) which was used to run the example of
Figure~\ref{fig:perturb}.
There also, as we have seen in section~\ref{pedagogical_example}, the
intermediate iterations are pertinent as they let us discover progressively the
smaller perturbations.

In conclusion, the {\overallbest} {\cuttingstrategy} seems to be well adapted
to image segmentation purposes, so we have used it in all the numerical results
below.

\subsection{Multiscalar segmentation}

For the general parameter estimation problem, multiscalar segmentation is a
natural approach when each component of the vector parameter has a different
physical interpretation.
This is not the case for image segmentation, where the three color components
have the same status.
But one can expect that allowing a different partition for each component
will allow to represent the image with a smaller number of regions for a same
level of details.
All tests in this section have been performed with the {\overallbest}
{\cuttingstrategy}.

\subsubsection{The {\bestcomponentonly} {\multiscalarstrategy}}

\begin{figure}[htb]
  \begin{center}
    (a)
    \includegraphics[width=\quartps\textwidth]{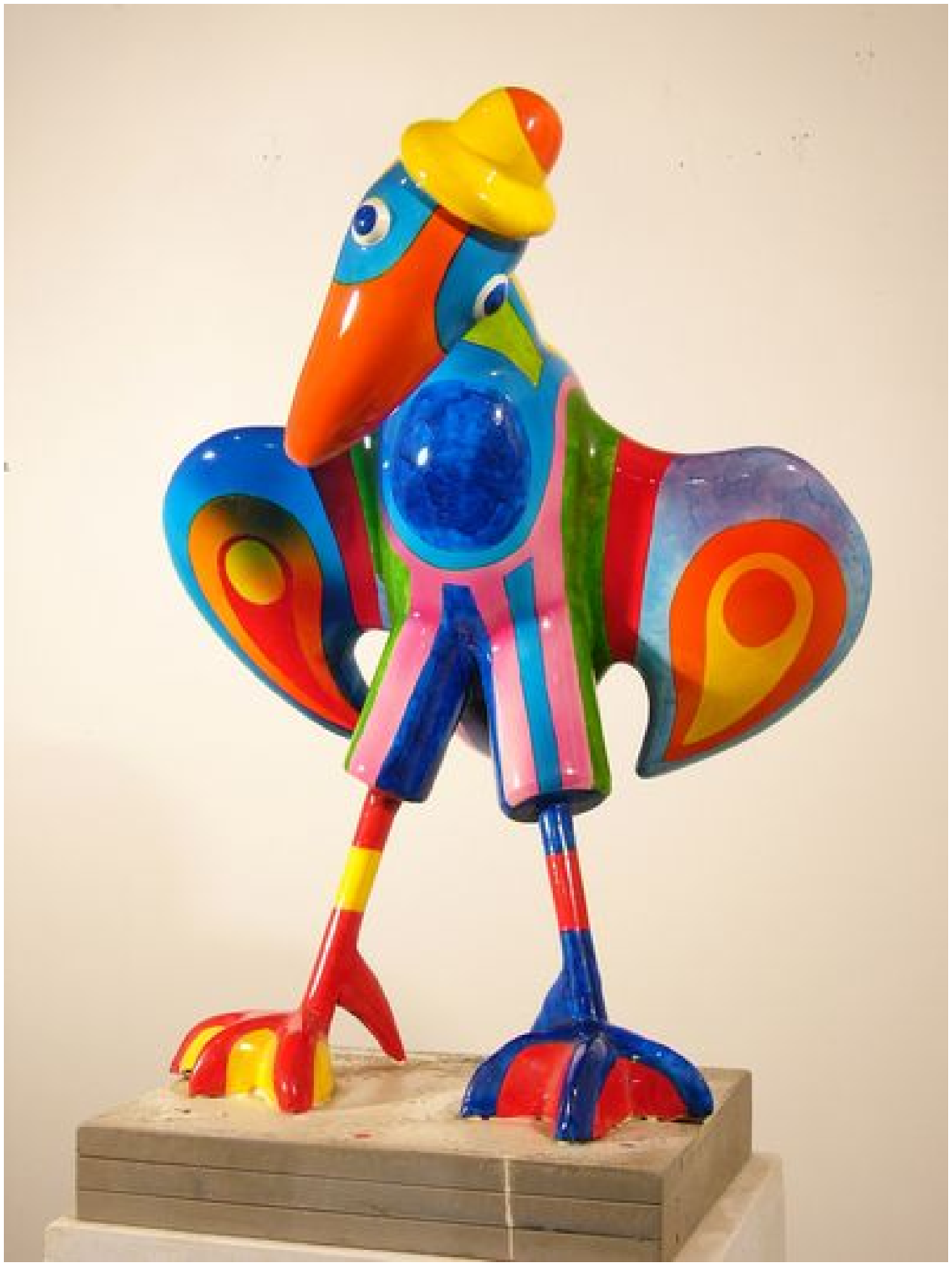}
    (b)
    \includegraphics[width=\quartps\textwidth]{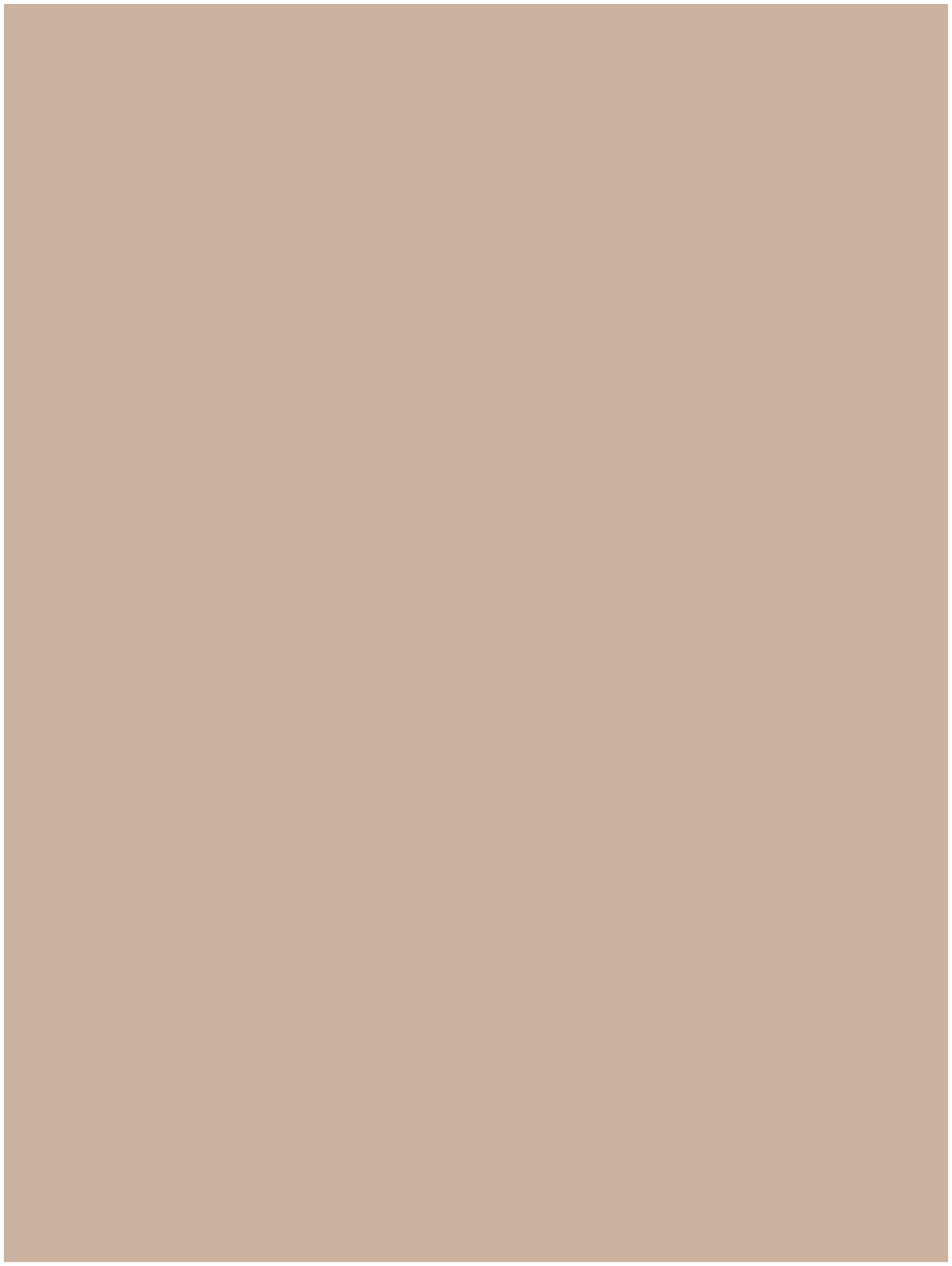}
    (c)
    \includegraphics[width=\quartps\textwidth]{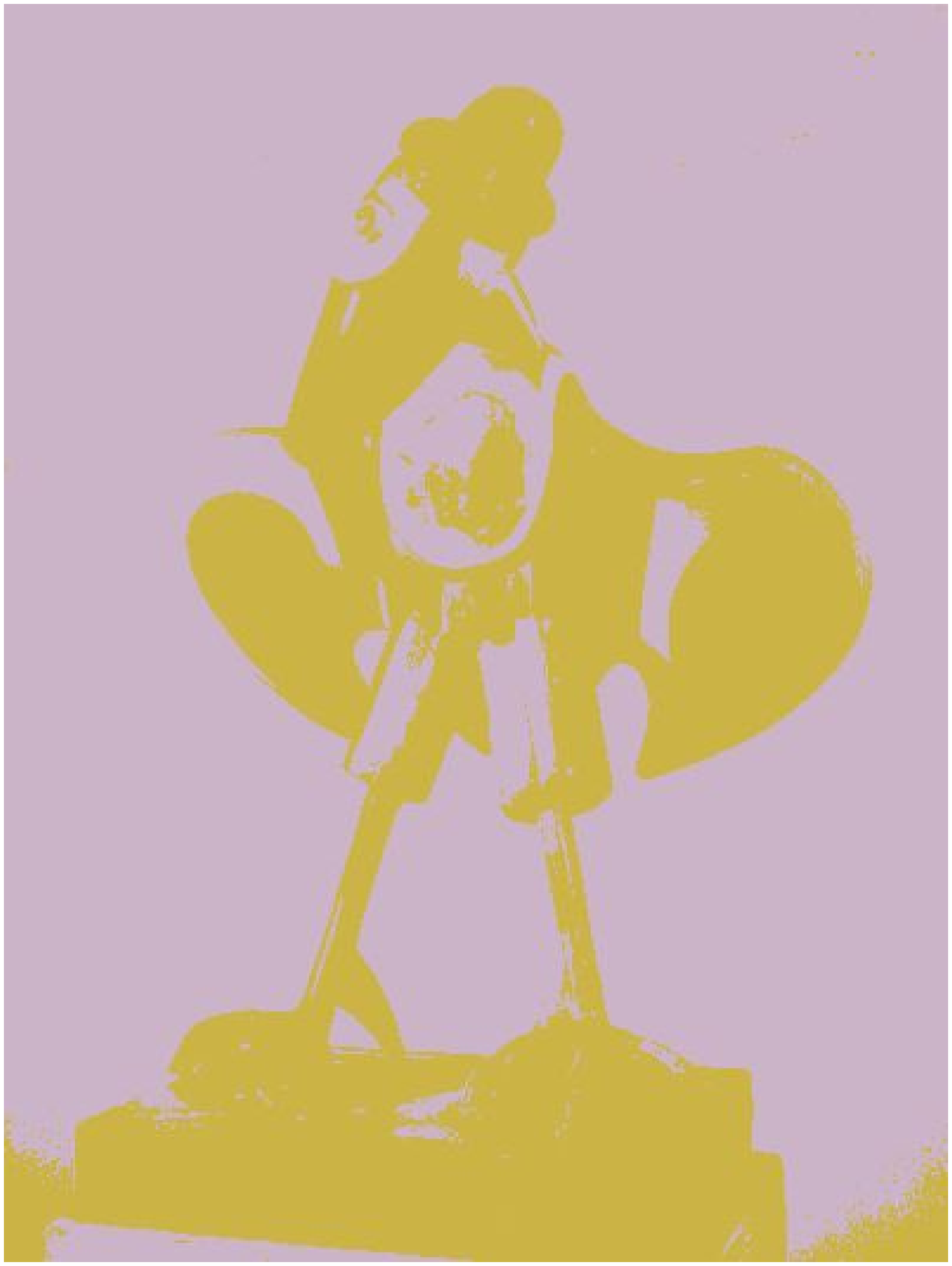} \\
    (d)
    \includegraphics[width=\quartps\textwidth]{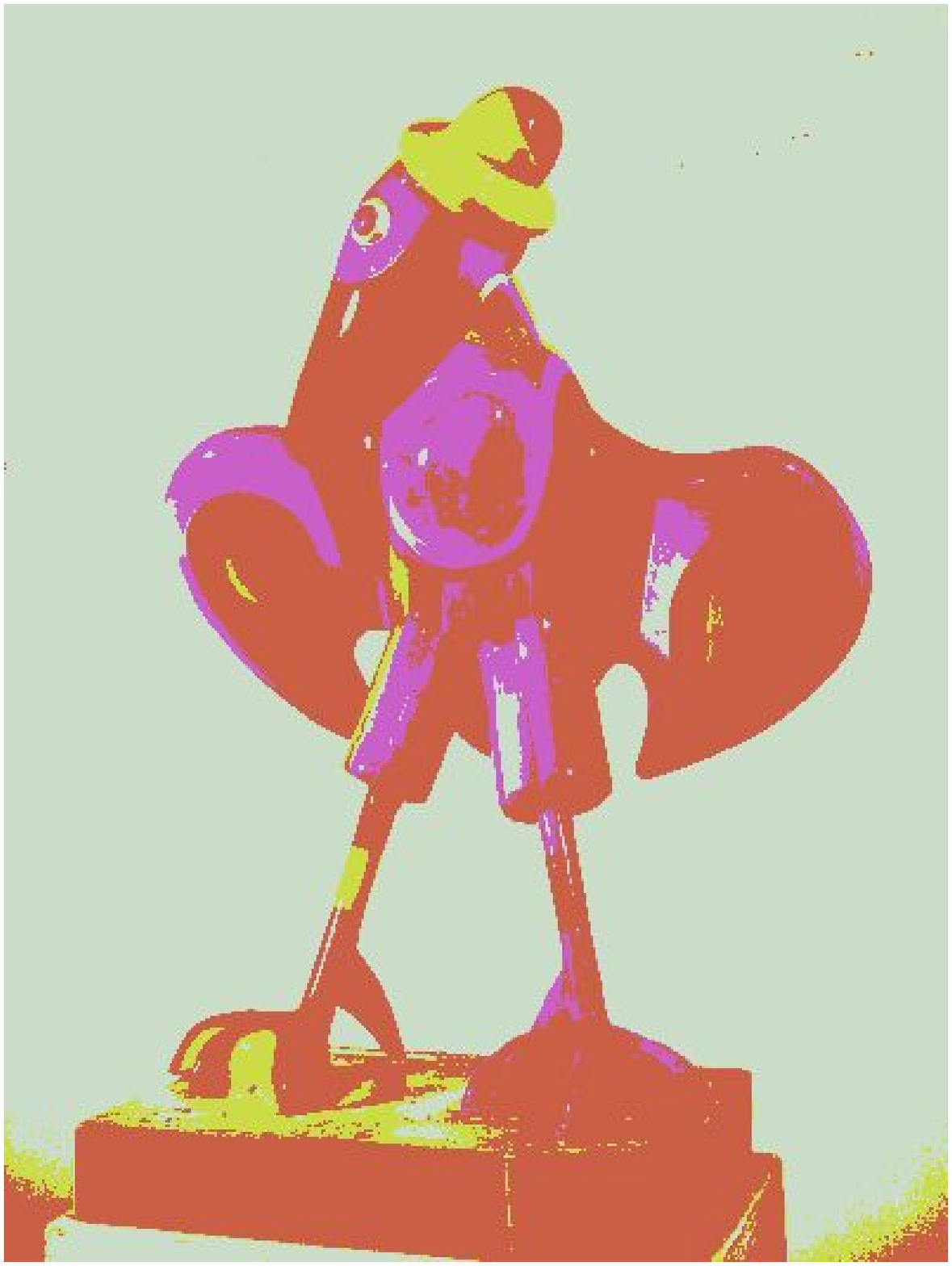}
    (e)
    \includegraphics[width=\quartps\textwidth]{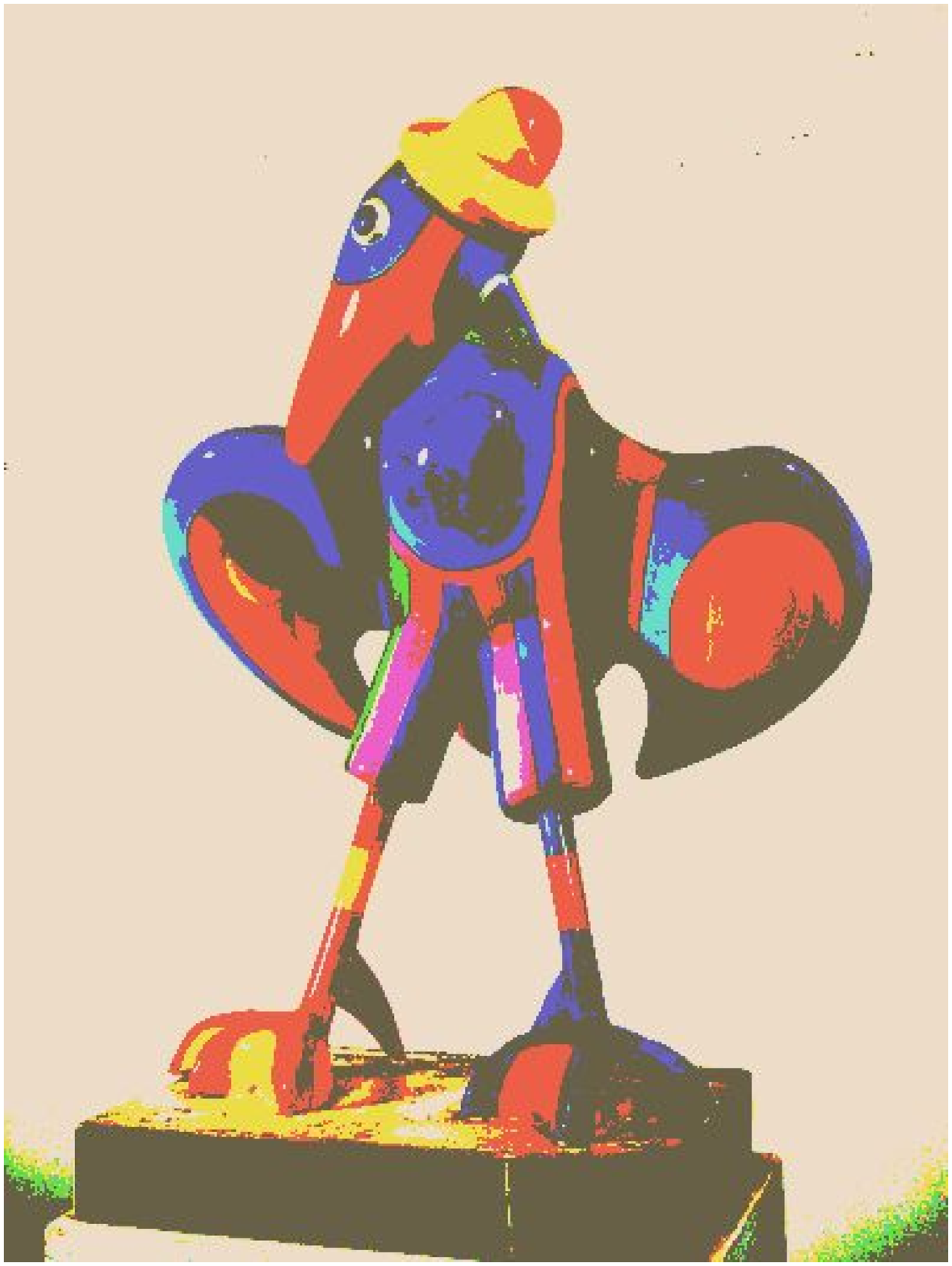}
    (f)
    \includegraphics[width=\quartps\textwidth]{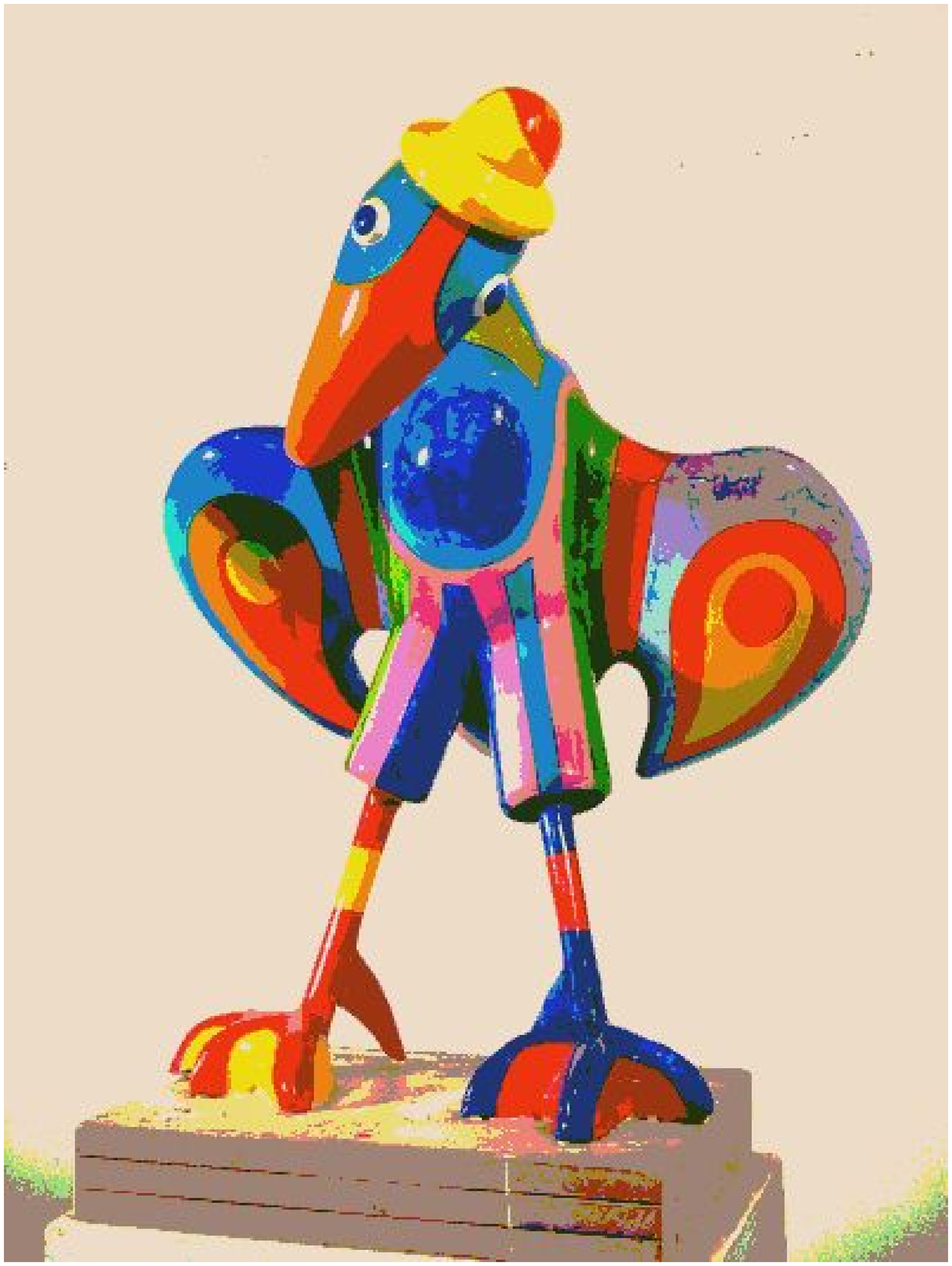}
    \caption{{\em Multiscalar segmentation} with the {\overallbest}
      {\cuttingstrategy} and the {\bestcomponentonly}
      {\multiscalarstrategy}.
      (a) is the image data to segment.
      Then, for each displayed image, we give the iteration number~$n$,
      the number of scalar regions~$n_{\rm sr}$, the resulting number of uniform
      color regions~$n_{\rm vr}$, and the percentage of explained
      data~$\tau$.\protect\\
      (b) $n=0$, $n_{\rm sr}=3$, $n_{\rm vr}=1$, $\tau=65.0\%$;\protect\\
      (c) $n=1$, $n_{\rm sr}=4$, $n_{\rm vr}=2$, $\tau=70.0\%$;\protect\\
      (d) $n=2$, $n_{\rm sr}=5$, $n_{\rm vr}=4$, $\tau=75.6\%$;\protect\\
      (e) $n=3$, $n_{\rm sr}=6$, $n_{\rm vr}=8$, $\tau=82.8\%$;\protect\\
      (f) $n=6$, $n_{\rm sr}=9$, $n_{\rm vr}=26$, $\tau=90.1\%$.}
    \label{fig:bestcomponentonly}
  \end{center}
\end{figure}

Here {\em one iteration} of the algorithm adds exactly
{\em one scalar degree of freedom} to the component which produces the
largest decrease of the objective function.
Hence, one can expect to obtain a good representation of the original image
with fewer component regions.
This is confirmed by the comparison of Figures~\ref{fig:bestcomponentonly}(f)
and~\ref{fig:overallbest}(g):
the {\bestcomponentonly} {\multiscalarstrategy} explains 90\%~of the data with
only 9~scalar regions (which generate 26~color regions), whereas the vector
segmentation algorithm requires 11~color regions, which correspond to
33~component regions.

\subsubsection{The {\bestcomponentforeach} {\multiscalarstrategy}}

\begin{figure}[htb]
  \begin{center}
    (a)
    \includegraphics[width=\quartps\textwidth]{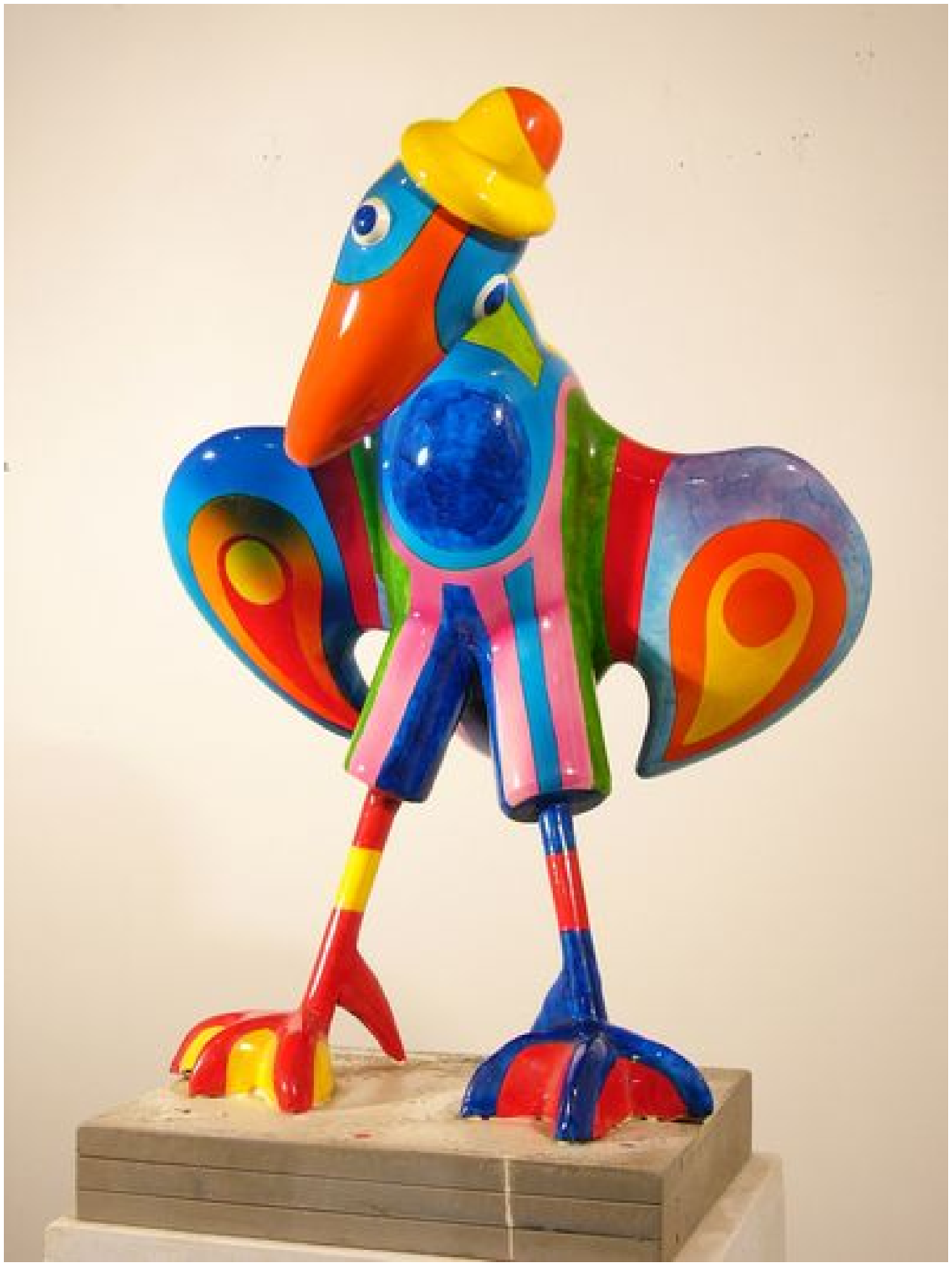}
    (b)
    \includegraphics[width=\quartps\textwidth]{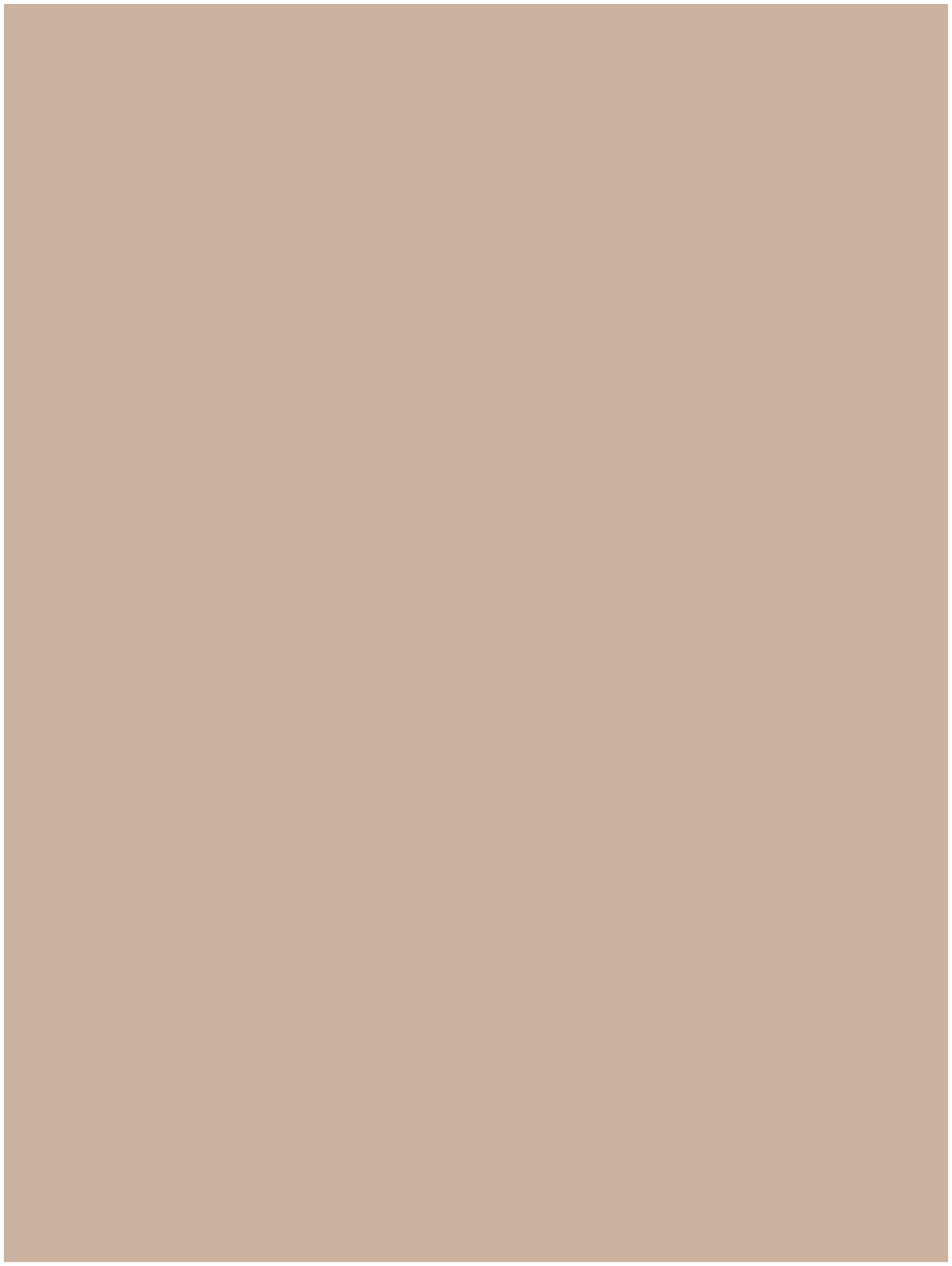}
    (c)
    \includegraphics[width=\quartps\textwidth]{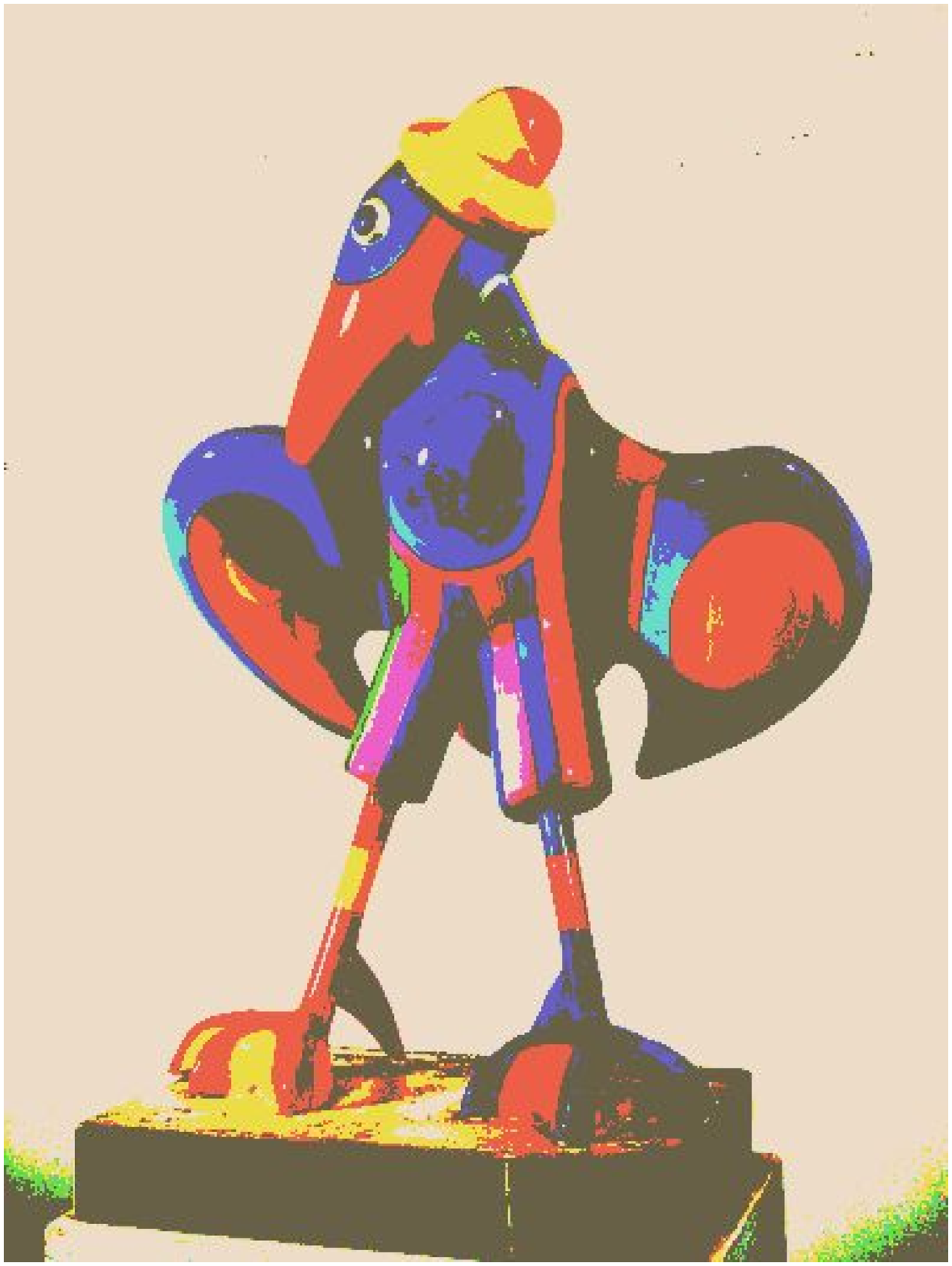} \\
    (d)
    \includegraphics[width=\quartps\textwidth]{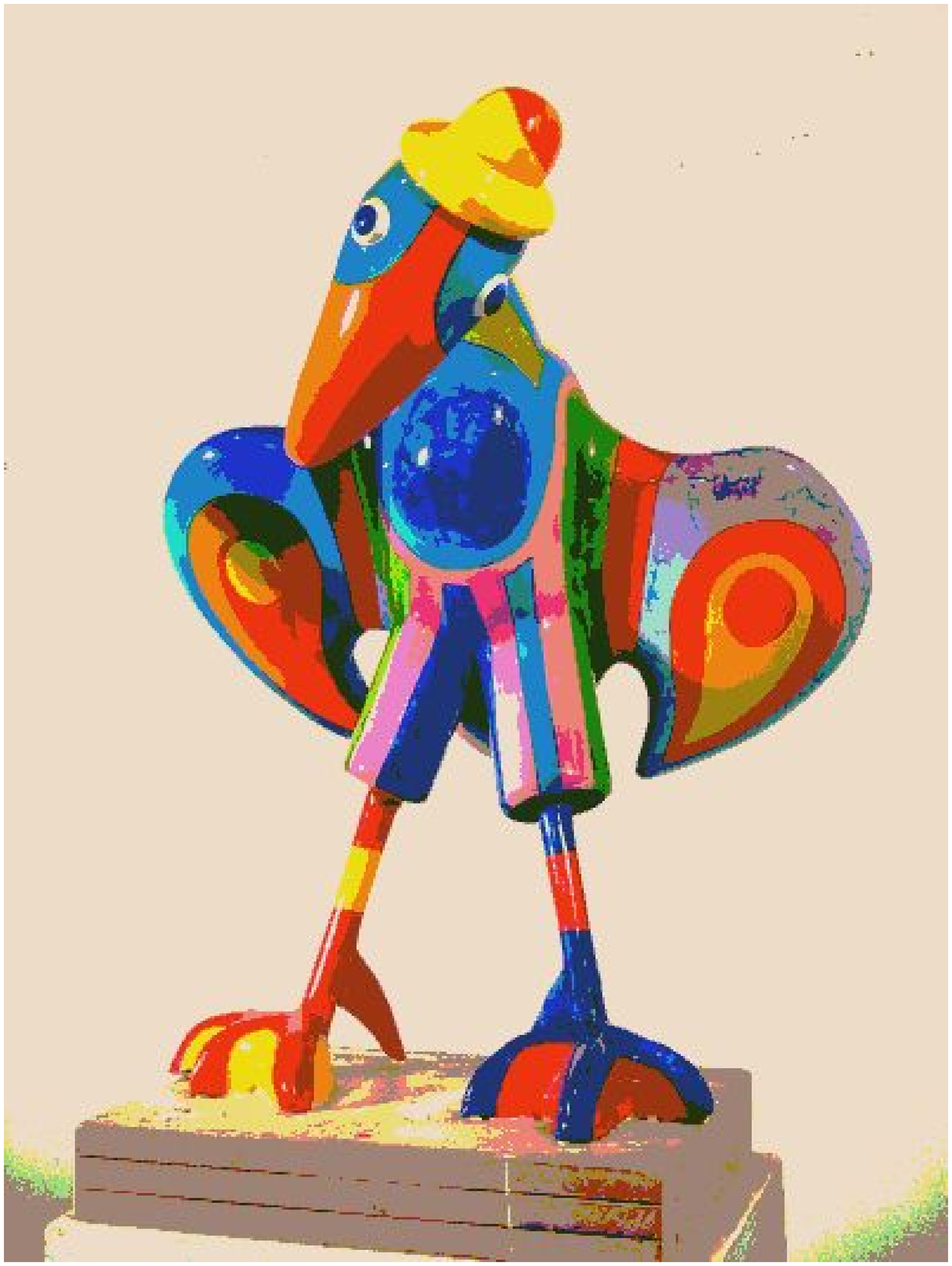}
    (e)
    \includegraphics[width=\quartps\textwidth]{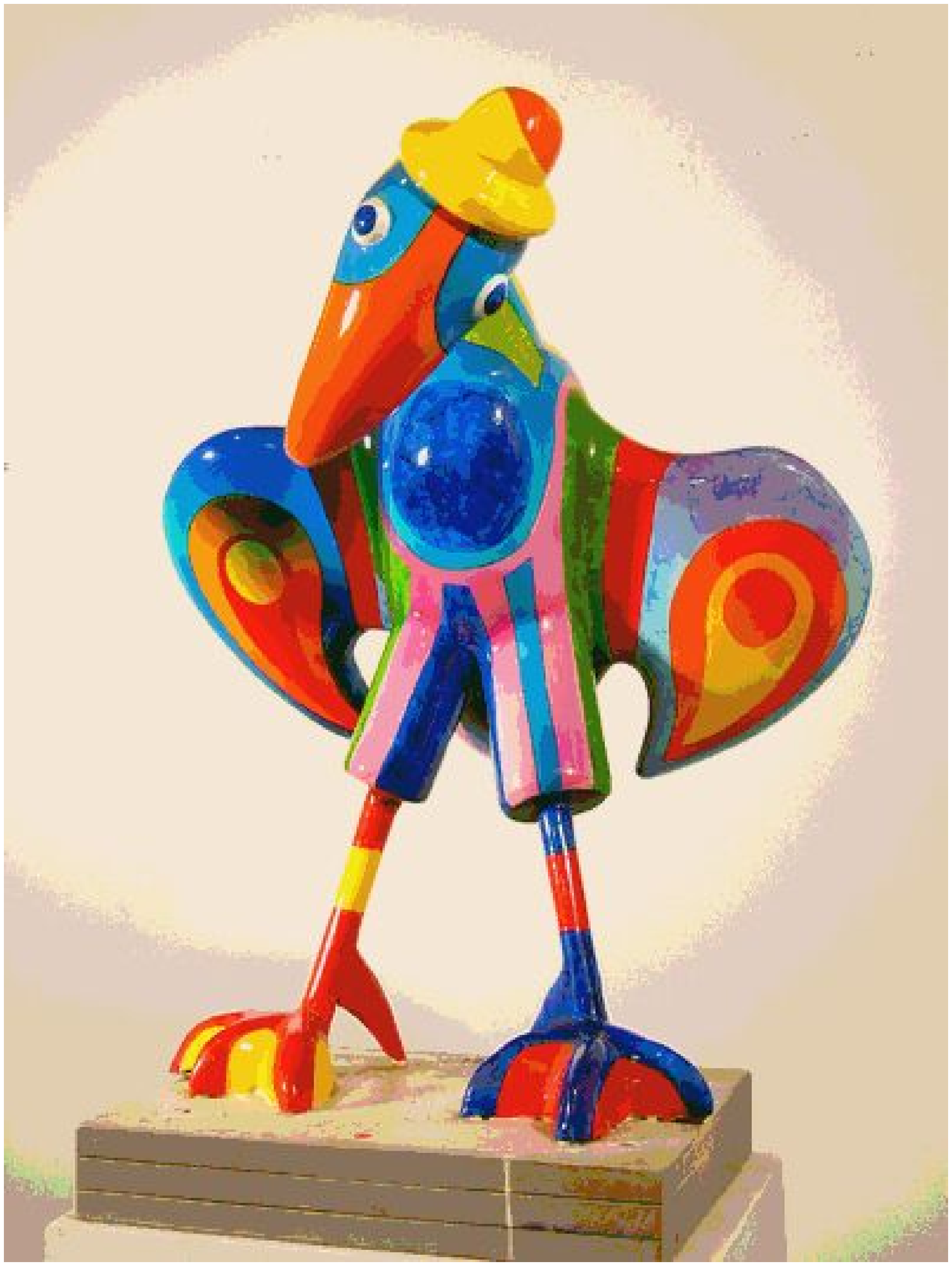}
    (f)
    \includegraphics[width=\quartps\textwidth]{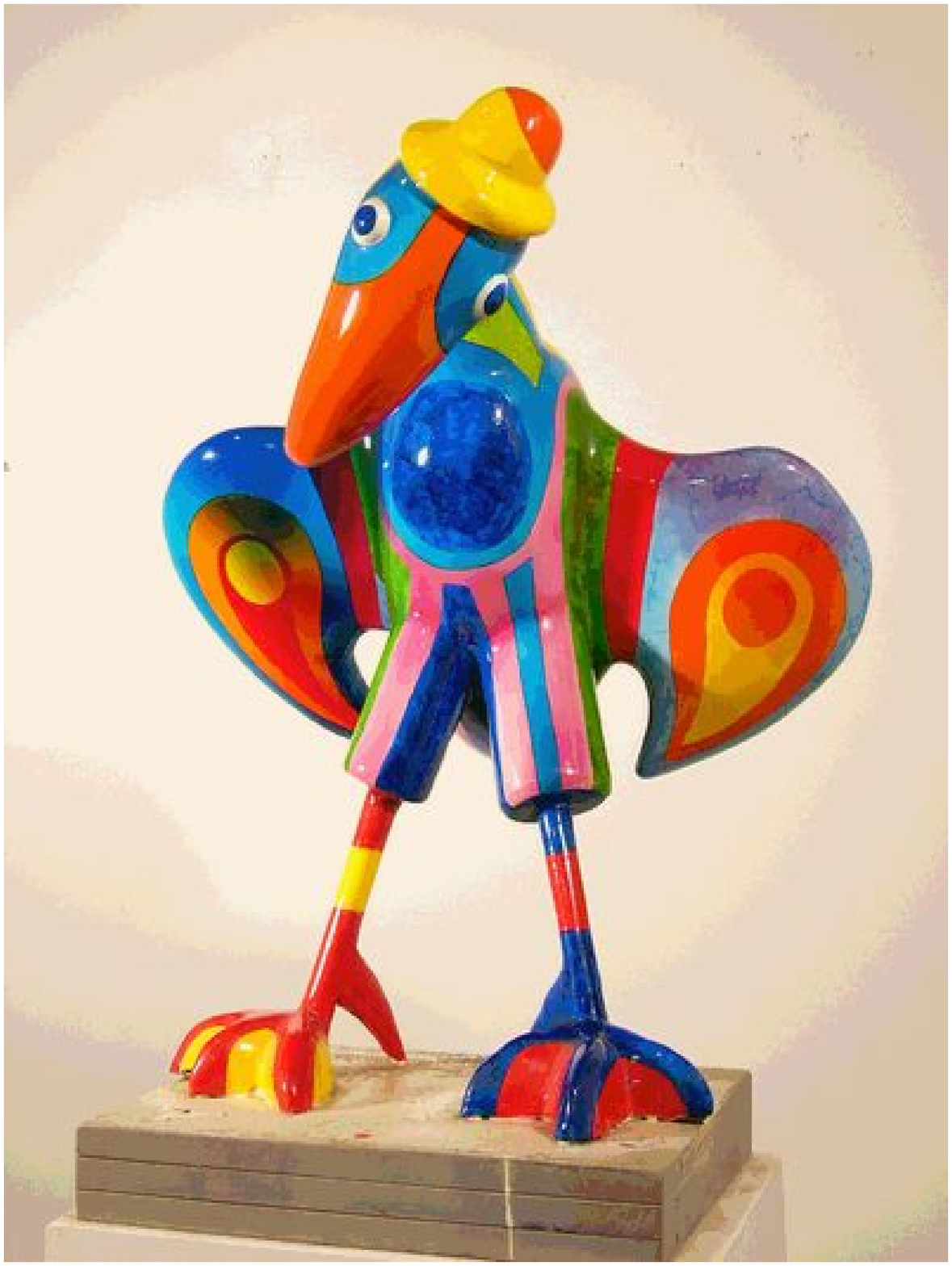}
    \caption{{\em Multiscalar segmentation} with the {\overallbest}
      {\cuttingstrategy} and the {\bestcomponentforeach}
      {\multiscalarstrategy}.
      (a) is the image data to segment.
      Then, for each displayed image, we give the iteration number~$n$,
      the number of scalar regions~$n_{\rm sr}$, the resulting number of vector
      regions~$n_{\rm vr}$, and the ratio of explained data~$\tau$.\protect\\
      (b) $n=0$, $n_{\rm sr}=3$, $n_{\rm vr}=1$, $\tau=65.0\%$;\protect\\
      (c) $n=1$, $n_{\rm sr}=6$, $n_{\rm vr}=8$, $\tau=82.8\%$;\protect\\
      (d) $n=2$, $n_{\rm sr}=9$, $n_{\rm vr}=26$, $\tau=90.1\%$;\protect\\
      (e) $n=5$, $n_{\rm sr}=18$, $n_{\rm vr}=153$, $\tau=94.7\%$;\protect\\
      (f) $n=10$, $n_{\rm sr}=33$, $n_{\rm vr}=4581$, $\tau=97.1\%$.}
    \label{fig:bestcomponentforeach}
  \end{center}
\end{figure}

The control on the number of degrees of freedom is less precise with this
algorithm, as {\em one iteration} adds exactly
{\em three scalar degrees of freedom} instead of one for the previous
algorithm.

Note that the images of Figures~\ref{fig:bestcomponentforeach}(d)
and~\ref{fig:bestcomponentonly}(f) coincide: this means simply that the
three first iterations---after the three initialization ones---of the
{\bestcomponentonly} algorithm have updated successively the~$R$, $G$
and~$B$ components.
But when the number of iterations increases, it can be expected that the
{\bestcomponentonly} algorithm will not cycle regularly through the $R$,
$G$ and~$B$ components  as the {\bestcomponentforeach} does, and hence produce
better images for a given number of component regions.

Because of the complete independence of the partitions for the~$R$, $G$ and~$B$
components, the number of color regions obtained by superimposition of the
component regions blows up for these two algorithms, see
Figure~\ref{fig:bestcomponentforeach}(f).
This shows that the parameterization of the image by  {\em component regions}
instead of {\em color regions} allows a drastic reduction of the number of
regions for a given image quality.
In other words, it is possible with these multiscalar algorithms to obtain a
segmented image with a large number of color regions, described by a much
smaller number of component regions.

\subsubsection{The {\combinebestcomponents} {\multiscalarstrategy}}

\begin{figure}[htb]
  \begin{center}
    (a)
    \includegraphics[width=\quartps\textwidth]{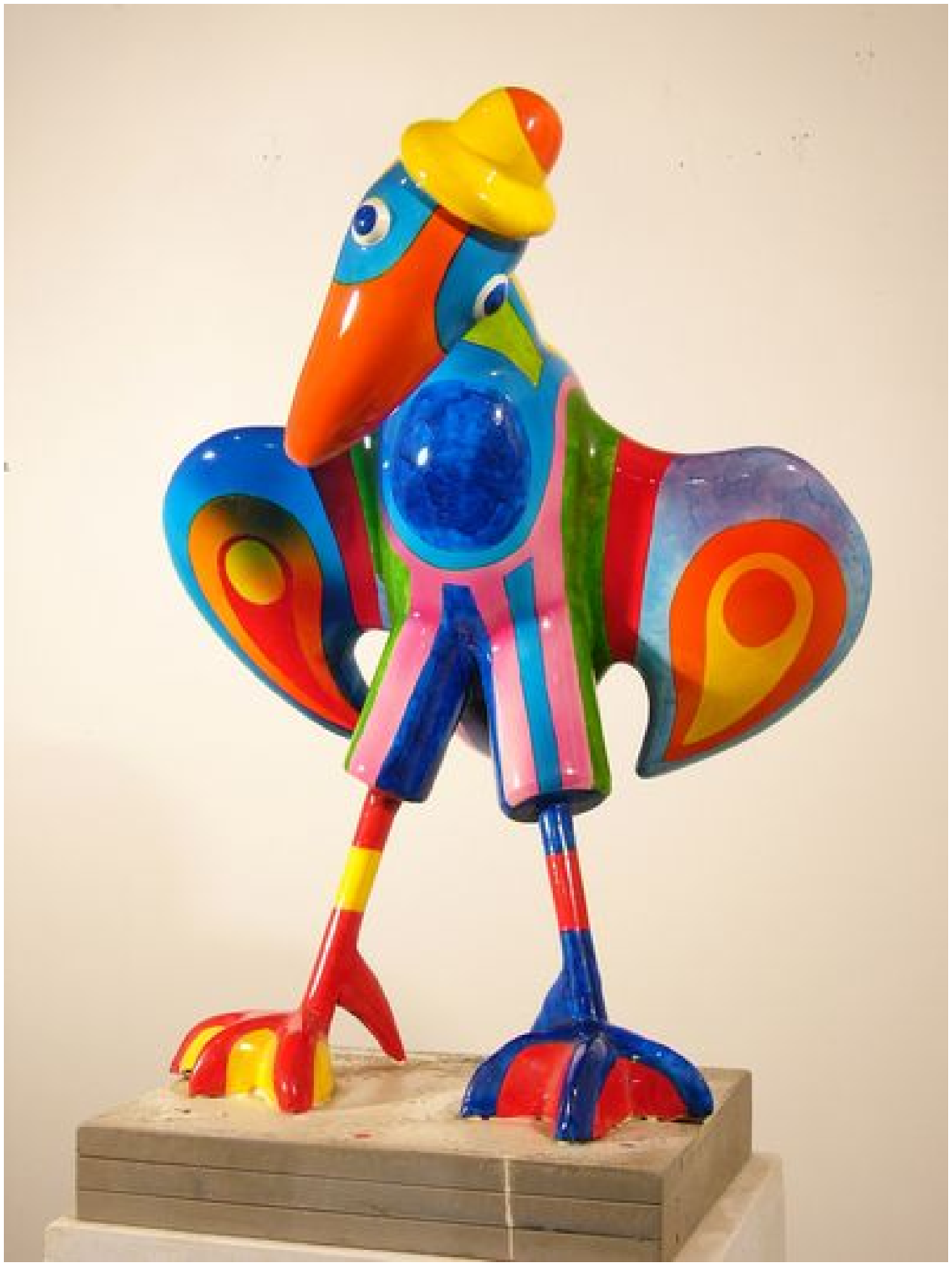}
    (b)
    \includegraphics[width=\quartps\textwidth]{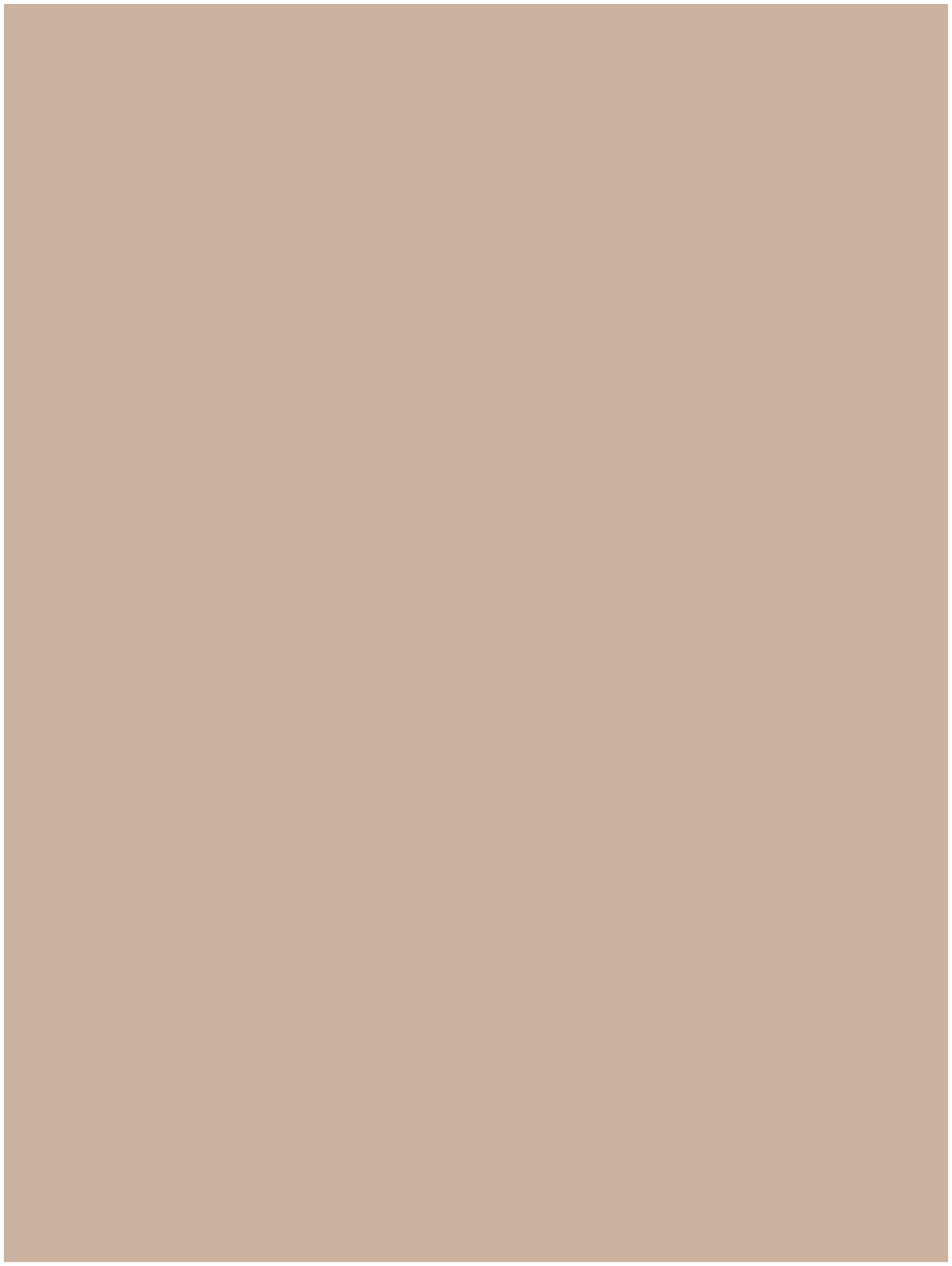}
    (c)
    \includegraphics[width=\quartps\textwidth]{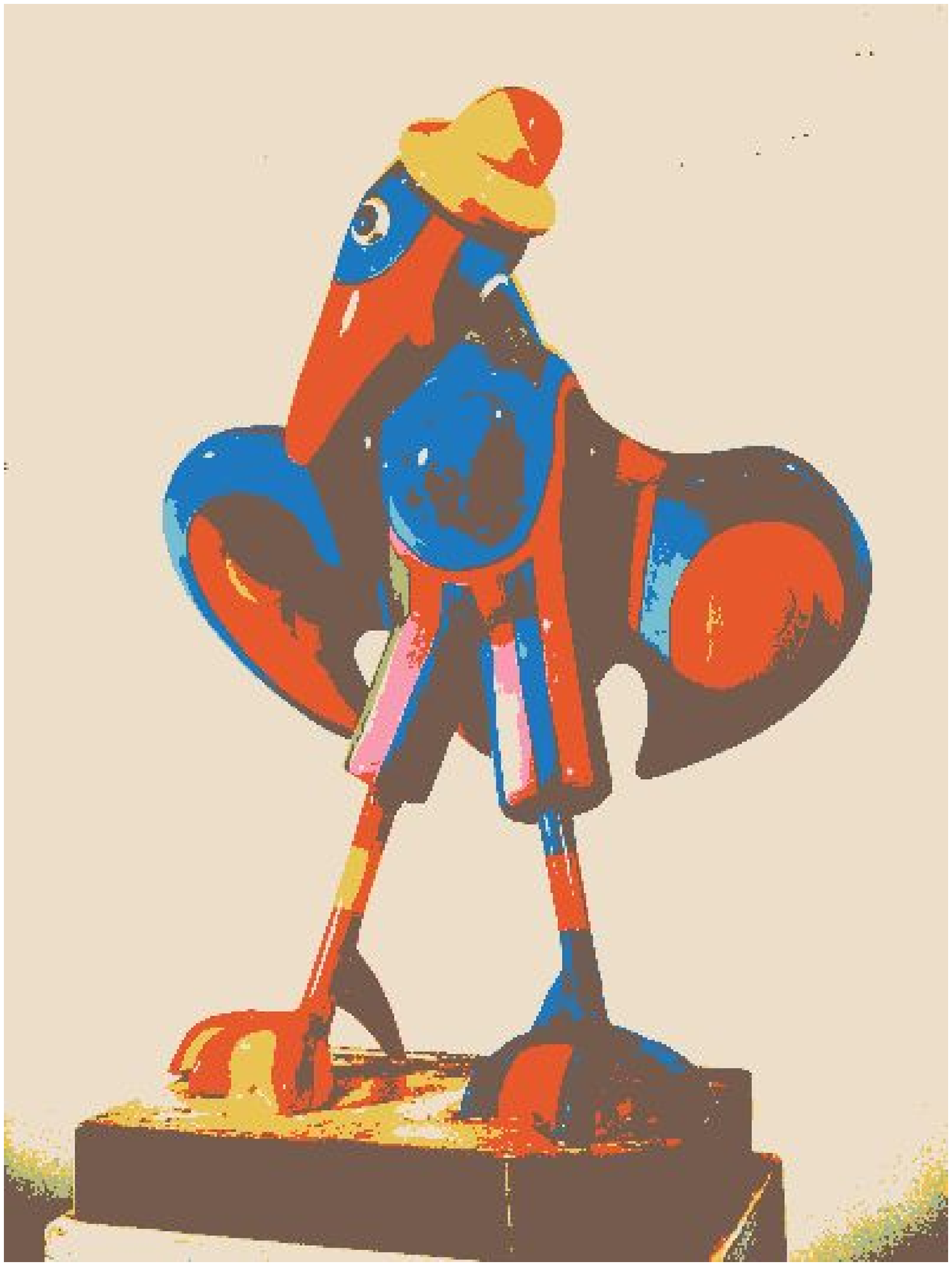} \\
    (d)
    \includegraphics[width=\quartps\textwidth]{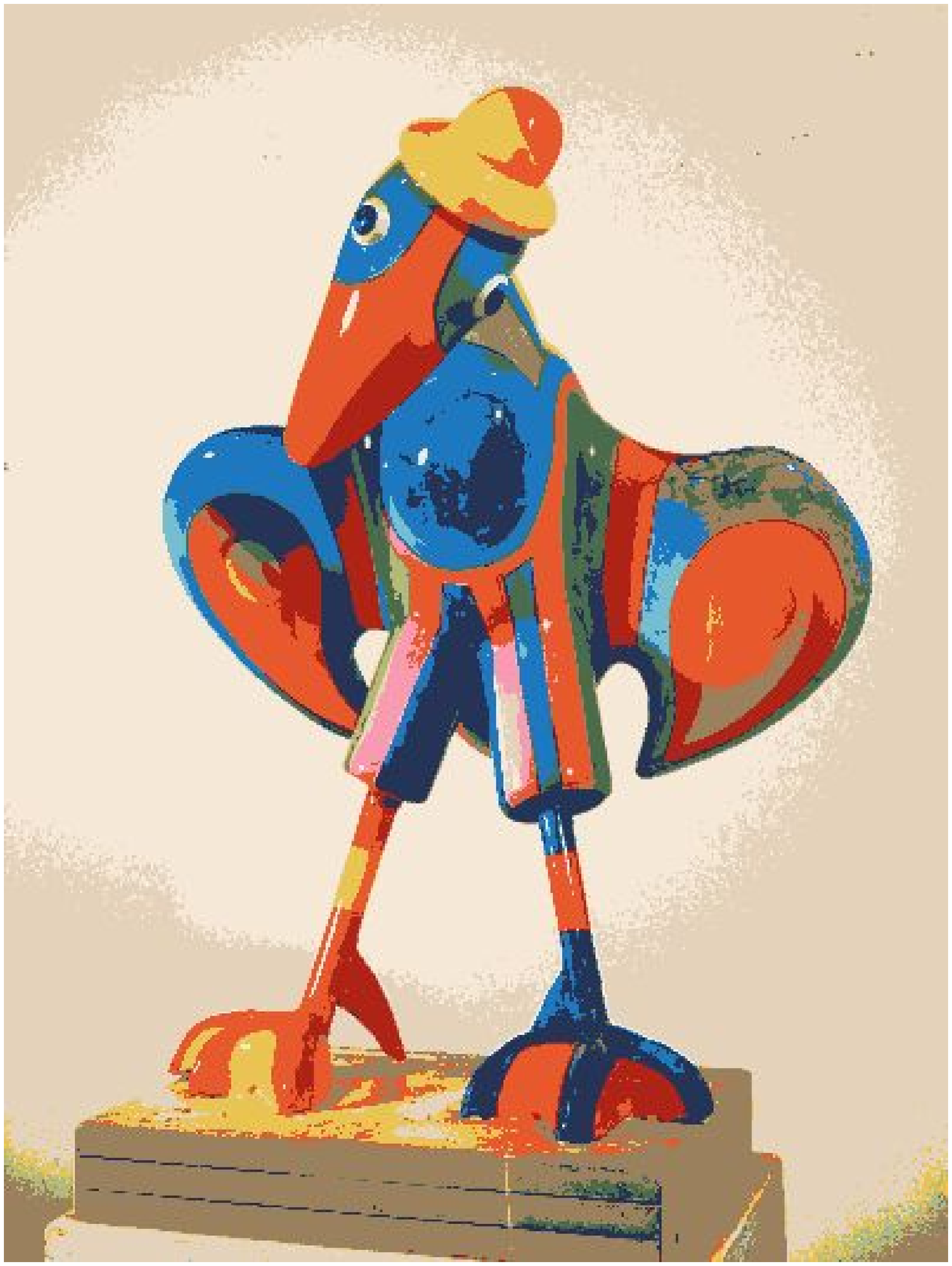}
    (e)
    \includegraphics[width=\quartps\textwidth]{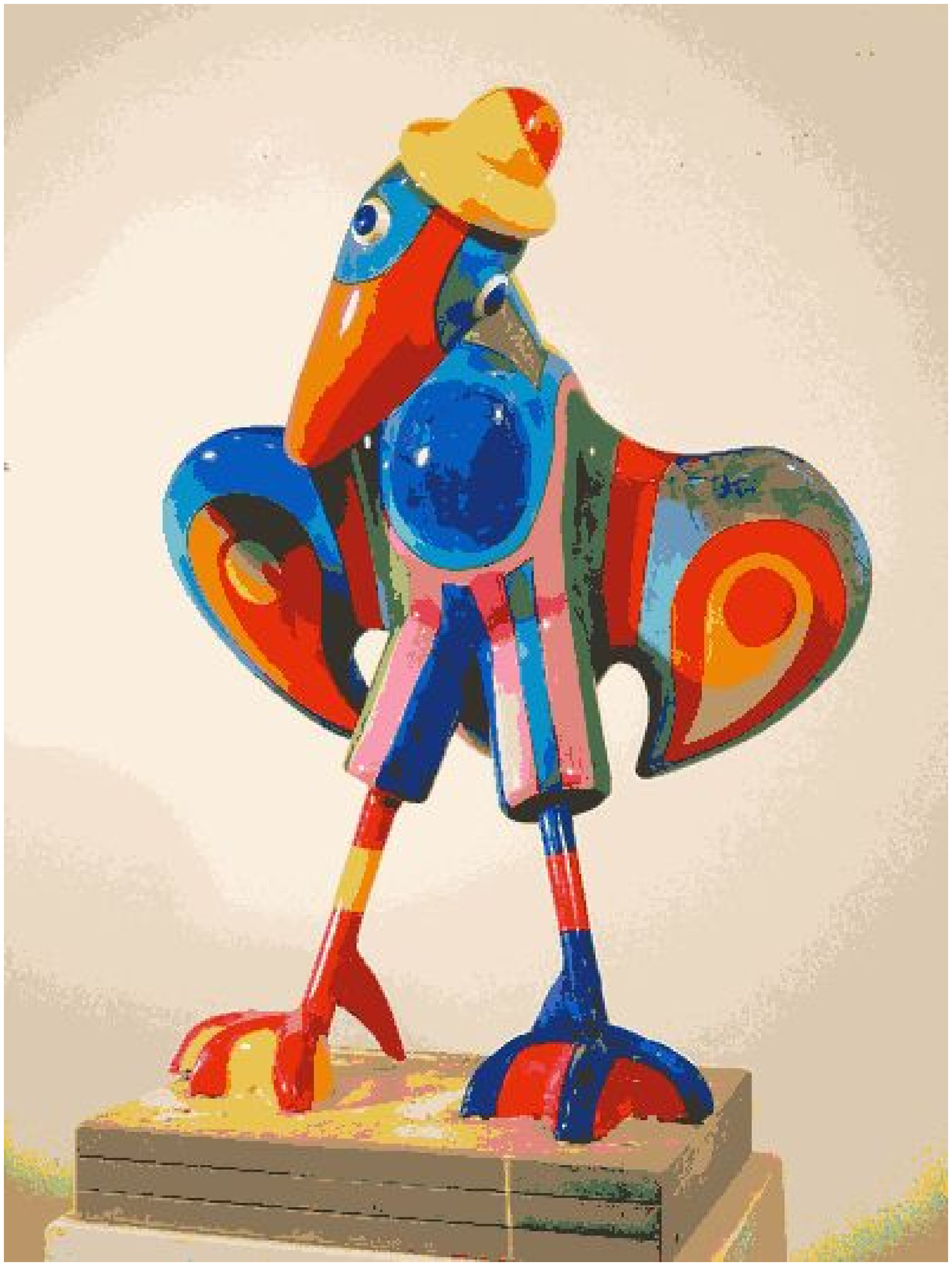}
    (f)
    \includegraphics[width=\quartps\textwidth]{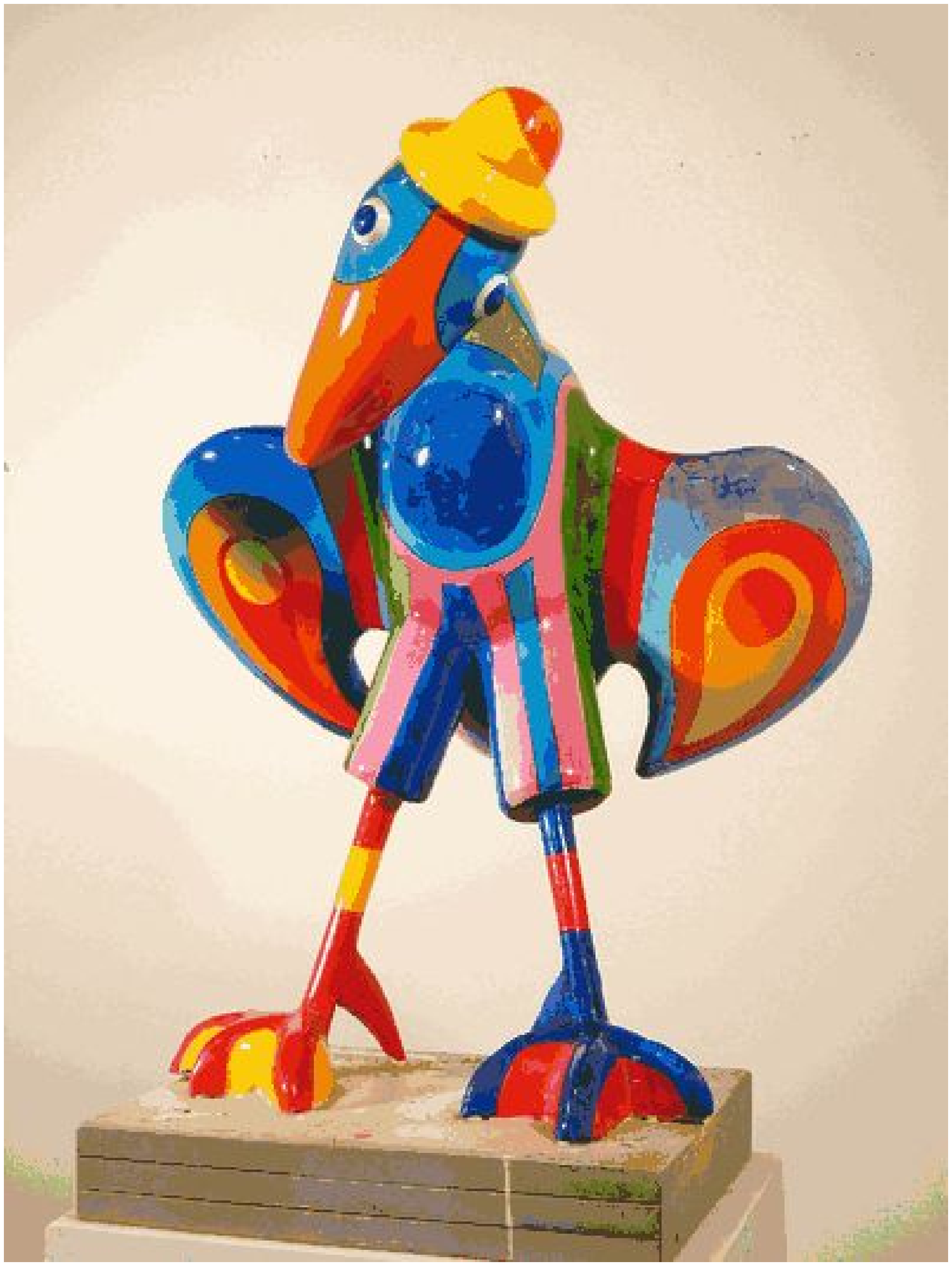}
    \caption{{\em Multiscalar segmentation} with the {\overallbest}
      {\cuttingstrategy} and the {\combinebestcomponents}
      {\multiscalarstrategy}.
      (a) is the image data to segment.
      Then, for each displayed image, we give the iteration number~$n$,
      the number of uniform color regions~$n_{\rm vr}$, and the percentage of
      explained data~$\tau$.\protect\\
      (b) $n=0$, $n_{\rm vr}=1$, $\tau=65.0\%$;\protect\\
      (c) $n=1$, $n_{\rm vr}=8$, $\tau=84.1\%$;\protect\\
      (d) $n=2$, $n_{\rm vr}=12$, $\tau=88.3\%$;\protect\\
      (e) $n=5$, $n_{\rm vr}=24$, $\tau=91.4\%$;\protect\\
      (f) $n=10$, $n_{\rm vr}=42$, $\tau=94.1\%$.}
    \label{fig:combinebestcomponents}
  \end{center}
\end{figure}

This {\multiscalarstrategy} produces in fact a vector segmentation, as we
have seen in the last item of section~\ref{sss:ms_strategy}.
Hence, all current component images are segmented on the same partition, so we
have omitted in Figure~\ref{fig:combinebestcomponents} the number~$n_{\rm sr}$ of
scalar regions, as it is given by $n_{\rm sr}=3n_{\rm vr}$.

Comparison of Figures~\ref{fig:combinebestcomponents}(d)
and~\ref{fig:combinebestcomponents}(e) with
Figure~\ref{fig:overallbest}(g) shows that, for a given level of
accuracy, the {\combinebestcomponents} {\multiscalarstrategy} tends to
produce vector segmented images with a larger number of regions than the
direct {\em vector segmentation algorithm}.

Also, the algorithm adds at each iteration a number of color regions between~3
and~7, i.e. between~9 and~21 scalar degrees of freedom, which makes it more
difficult to control of the number of color regions.

So in a whole, the {\combinebestcomponents} {\multiscalarstrategy} seems
less performing for image {\em vector segmentation} than the direct vector
segmentation algorithm~\ref{sss:algo_vector_segmentation}.


\section*{Conclusion}

Once recognized as a parameter estimation problem, for which the forward
modeling map to invert is the identity, image segmentation may benefit from
the theory of optimal control.

We have studied the specificities of this favorable case in the context of
adaptive parameterization and we have proposed an
{\em optimal adaptive (vector) segmentation algorithm}.
This algorithm provides a finite sequence of segmented images with a number of
{\em color regions} increasing by one at each iteration.
Thanks to its adaptive character, the algorithm extracts the significant
regions from the very beginning, and all the intermediate steps of the process
are meaningful.
The optimal character indicates that the {\em most significant} regions are
extracted first, and that the sequence minimizes the error with respect to
the original image and terminates with a segmented image containing as many
regions as distinct colors in the original image corresponding to the null
error.
This algorithm is well-suited to finely control the desired number of color
regions, since it is equipped with a meaningful criterion to stop the
refinement process.

The study has also exhibited a class of {\em multiscalar algorithms} that
provide a distinct segmentation for each color component.
Here the number of {\em component regions} increases by one or three at each
iteration, but the number of resulting {\em color regions} blows up quickly,
thus explaining a large percentage of the data.
Hence, this algorithm also shines when a desired percentage of the data has to
be explained with a number of (component) regions as small as possible.

\section*{Acknowledgements}

This research was partially carried out as part of the DGRSRT/INRIA STIC
project ``Identification de param\`etres en milieu poreux : analyse
math\'ematique et \'etude num\'erique'', and of the ENIT/INRIA Associated Team
MODESS.
The authors hereby acknowledge the provided support.

The authors would also like to thank Les Popille for their kind permission to
use a reproduction of one of their sculptures.

\bibliography{biblio}
\bibliographystyle{plain}

\end{document}